\documentclass[brochure,12pt]{bourbaki}% Pour un exposŽ en franais
\usepackage[matrix,arrow]{xy}
\usepackage{amssymb,amsfonts,amsmath,footnote}
\usepackage[francais]{babel}
\usepackage{graphics}
\usepackage{graphicx}
\usepackage{psfrag}
\usepackage{PSFigure,psfrag}
\usepackage{etex}
\usepackage{color}
\usepackage{tikz-cd}
\newcommand{\be}{\begin{equation}}
\newcommand{\ee}{\end{equation}}
\newcommand{\R}{\mathbb{R}}
\newcommand{\N}{\mathbb{N}}

\newcommand{\Z}{\mathbb{Z}}

\def\thesection{\Roman{section}}
\def\theequation{\thesection.\arabic{equation}}

\def\ti{\tilde}
\def\lf{\left}
\def\rg{\right}

\def\al{\alpha}
\def\la{\lambda}

\def\ep{\varepsilon}
\def\ds{\displaystyle}
\def\ov{\overline}
\def\Om{\Omega}

\def\p{\partial}

\date{Juin 2019}
\bbkannee{71\`eme ann\'ee, 2018-2019}
\bbknumero{1163}
\title{Infinit\'e d'hypersurfaces minimales en basses dimensions}
\subtitle{d'apr\`es F.C. Marques, A.A. Neves  et A. Song}
\author{Tristan RIVI\`ERE}
\address{ETH-Zentrum\\
Forschungsinstitut f\"ur Mathematik\\
CH--8092 Z\"urich, Switzerland}
\email{tristan.riviere@math.ethz.ch}

\begin{document}
\maketitle

\noindent{\bf INTRODUCTION}

\bigskip

La recherche de g\'eod\'esiques ferm\'ees, c'est \`a dire des cercles points critiques de la fonctionnelle de longueur, sur une surface ou une vari\'et\'e de dimension quelconque est au moins aussi ancienne que le calcul de Euler et Lagrange entre 1744 et 1760. Cette recherche trouve une g\'en\'eralisation naturelle en dimension plus grande dans la notion de surfaces minimales qui sont les points critiques du volume parmi les sous vari\'et\'es de dimension et de topologie donn\'ee. Les probl\'ematiques consistant \`a mettre en \'evidence des ``sous vari\'et\'es optimales''   \`a un espace donn\'e se sont av\'er\'ees \^etre des moteurs f\'econds dans le d\'eveloppement de pans entiers des math\'ematiques comme le calcul des variation, les \'equations aux d\'eriv\'ees partielles, la g\'eom\'etrie  et la topologie diff\'erentielle, les syst\`emes dynamiques...etc. 

Un r\'esultat  classique de Marston Morse (voir par exemple \cite{Kl}) affirme que sur certains ellipsoides de ${\R}^3$ il existe exactement 3 g\'eod\'esiques ferm\'ees et simples (plong\'ees).
Le but de cet expos\'e est d'\'etablir qu'un tel r\'esultat de rigidit\'e ne s'\'etend pas aux dimensions sup\'erieures et que dans toute vari\'et\'e riemannienne lisse, compacte et sans bord de dimension inf\'erieure ou \'egale \`a 7, il existe une infinit\'e d'hypersurfaces minimales plong\'ees distinctes\footnote{La perte de cette rigidit\'e en dimension sup\'erieure doit \^etre relativis\'ee : si on se restreint aux surfaces minimales de topologie donn\'ee alors  on a par exemple le r\'esultat suivant : il existe des m\'etriques voisines de la m\'etrique standard sur $S^3$ qui admettent exactement 4 sph\`ere minimales plong\'ees (\cite{Whi}). Un r\'esulat de J\"urgen Jost  dit par ailleurs que toute sph\`ere tridimensionnelle admet au moins 4 sph\`eres minimales plong\'ees  \cite{Jos}}. Ce th\'eor\`eme dont nous allons d\'ecrire les origines et la d\'emonstration r\'ecente  par Antoine Song est l'aboutissement d'efforts qui nous renvoient  aux origines du calcul des variations et des probl\`emes de {\it minmax} avec  les travaux de Georges Birkoff, Lazar Lyusternik, Lev Shnirelmann, Abram Fet, Marsden Morse, Richard Palais, Stephen Smale...etc,  tout comme \`a la fusion de la th\'eorie de la mesure g\'eom\'etrique\footnote{ La th\'eorie de la mesure g\'eom\'etrique a \'et\'e motiv\'ee \`a l'origine par la r\'esolution du probl\`eme de Plateau dans sa plus grande g\'en\'eralit\'e.} avec cette th\'eorie de minmax sous l'impulsion de Fred Almgren et Jon Pitts.

\section{La construction de g\'eod\'esiques par minmax}

Les approches pour construire des g\'eod\'esiques ferm\'ees sont nombreuses et tr\`es riches. Une pr\'esentation exhaustive de celles-ci et des
th\'eories math\'ematiques qui en d\'ecoulent va bien au del\`a de cet expos\'e. Nous nous contentons de mentionner l'origine des id\'ees qui ont directement influenc\'e la preuve qui nous int\'eresse aujourd'hui.

Suite au travaux de Jacques Hadamard et d'Henri Poincar\'e sur le sujet, Georges Birkhoff a d\'evelop\'e \`a la fois des approches de syst\`emes dynamiques (flot g\'eod\'esique) mais aussi variationnelles afin de construire des g\'eod\'esiques ferm\'ees dans des vari\'et\'es simplement connexes (\cite{Bir}). Ces derni\`eres
peuvent \^etre revisit\'ees dans l'esprit de la th\'eorie de Morse sur le sujet ainsi qu'au moyen de la th\'eorie des d\'eformations en dimension infinie
de Richard  Palais et Stephen Smale, dont elles sont \`a l'origine. 

Pr\'ecisemment, \'etant donn\'ee une sous-vari\'et\'e riemannienne ferm\'ee (compacte sans bord) et connexe $(M^m,g)$ de ${\R}^Q$, dans le cas o\`u $\pi_1(M^m)\ne 0$ l'existence d'une g\'eod\'esique ferm\'ee s'obtient par simple minimisation de l'\'energie de Dirichlet
\[
E(u):=\int_0^{2\pi} \lf|\frac{\p u}{\p \theta}\rg|^2\ d\theta
\]
dans le sous ensemble de la vari\'et\'e de Hilbert 
\[
{\mathcal M}:=W^{1,2}(S^1,M^m):=\lf\{u\in W^{1,2}(S^1,{\R}^Q)\ ;\ u(\theta)\in M^m\ ,\  \forall \theta\in S^1\rg\}
\]
des \'el\'ements $u$ r\'ealisant une classe de conjugaison donn\'ee du $\pi_1$ de $M^m$. L'injection compacte  $W^{1,2}(S^1,M^m)\hookrightarrow C^0(S^1,M^m)$ garantit
l'existence d'un minimum dans chacune de ses classes qui s'av\`ere \^etre une g\'eod\'esique ferm\'ee param\'etris\'ee \`a vitesse constante.

Lorsque la vari\'et\'e $M^m$ est simplement connexe l'obtention d'une g\'eod\'esique ferm\'ee s'obtient par un argument variationnel plus \'elabor\'e qu'une simple minimization : un argument dit de {\bf minmax}.

Un r\'esultat classique de la th\'eorie de l'homotopie des vari\'et\'es compactes donne l'existence de $k\in \{2,\cdots, m\}$ tel que $\pi_k(M^m)\ne0$ mais $\pi_l(M^m)=0$ pour $l\in \{1\cdots k-1\}$. On d\'emontre sans trop d'efforts que la vari\'et\'e de Hilbert ${\mathcal M}$ est homotopiquement \'equivalente \`a l'espace des chemins ferm\'es $C^0(S^1,M^m)$. On d\'enote $\Om_p(M^m)$ l'espace des lacets de $M^m$, c'est-\`a-dire l'ensemble des \'el\'ements $u$ dans  $C^0(S^1,M^m)$ tels que $u(0)=p$ o\`u $p$ est un point fix\'e de $M^m$. On observe que l'application ``\'evaluation'' $ev$ de  $C^0(S^1,M^m)$ dans $M^m$ qui \`a $u$ associe $u(0)=u(2\pi)$ est  une fibration de Serre de fibre l'espace des lacets $\Omega(M^m)$.  Cette projection poss\`ede un inverse $i$ qui \`a un \'el\'ement quelconque de $M$ associe le chemin constant \'egal \`a cet \'el\'ement. La suite exacte de fibration
\[
\cdots \longrightarrow \pi_n(\Om_p(M^m)) \longrightarrow 
\begin{tikzcd}\pi_n(C^0(S^1,M^m)) \arrow[r, "\mathrm{ev}_{\ast}"'] &   \ar[l,bend left=80,looseness=1.8, "\iota_{\ast}"] \pi_{n}(M^m)
\end{tikzcd}\longrightarrow \pi_{n-1}(\Om_p(M^m))\longrightarrow\cdots
\]
se coupe en morceaux (i.e. Im $ev_\ast=\pi_{n}(M^m)$, en effet tout \'element de $\pi_{n}(M^m)$ est envoy\'e par $i_\ast$ en un \'el\'ement de $\pi_n(C^0(S^1,M^m))$ et l'on a clairement $ev_\ast\circ i_\ast=id$ ) et donne l'identit\'e
\[
\pi_n(C^0(S^1,M^m))=\pi_n(\Om_p(M^m))\oplus \pi_{n}(M^m)
\]
On a tautologiquement $\pi_n(\Om_p(M^m))=\pi_{n+1}(M^m)$ ce qui donne finalement
\[
\pi_n({\mathcal M})=\pi_{n+1}(M)\oplus \pi_{n}(M^m)
\]
 On a donc en particulier pour l'entier $k$ d\'efini plus haut : $\pi_{k-1}({\mathcal M})\ne 0$.\footnote{Par exemple, si $M^m=S^2$ on a $$\pi_1(W^{1,2}(S^1,S^2))={\Z}$$ et un g\'en\'erateur de ce premier groupe 
 d'homotopie est donn\'e par ce que l'on appelle un {\bf balayage de Birkhoff} : la famille continue de chemins sur $S^2$ obtenue en prenant l'intersection de cette sph\`ere avec tous les
 plans horizontaux  de ${\R}^3$``balayant'' de haut en bas $S^2$ exactement une fois, ensuite la boucle est referm\'ee au moyen dune famille continue de chemins constant reliant le p\^ole sud au p\^ole nord.}
 
 On introduit alors le {\bf probl\`eme de minmax} suivant :  existe -t-il un point critique non trivial (non \'egal \`a une constante) de $E$ qui r\'ealise le nombre donn\'e par
 \be
 \label{1}
 W:=\inf_{\lf\{u(s,\theta)\in C^0(S^{k-1},{\mathcal M})\ ; \ [u]\ne 0\rg\}}\quad \max_{s\in S^{k-1}}E(u(s,\cdot))\ \quad ?
 \ee
 o\`u $[u]$ d\'esigne la classe d'homotopie de $u$. Le nombre $W$ est appel\'e la {\bf largeur} du probl\`eme de minmax. La r\'eponse \`a cette question est positive et sa d\'emonstration donne l'existence d'une g\'eod\'esique ferm\'ee
 non triviale dans toute vari\'et\'e riemanniene ferm\'ee r\'eguli\`ere\footnote{Toute vari\'et\'e ferm\'ee et lisse se plonge isom\'etriquement dans un espace euclidien grace \`a un th\'eor\`eme classique d\^u \`a John Nash.}.
 
 Il faut d'abord s'assurer que $W>0$. Si tel n'\'etait pas le cas il existerait alors pour tout $\ep>0$ un \'el\'ement $u_\ep\in C^0(S^{k-1},{\mathcal M})$ tel que $[u_\ep]\ne 0$ et
 \[
 \max_{s\in S^{k-1}}E(u_\ep(s,\cdot))<\ep
 \] 
 En choisissant $\ep$ suffisament petit (inf\'erieur au rayon de convexit\'e de $M^m$) on aurait ainsi le fait que pour tout $t\in S^{k-1}$ le chemin $u_\ep(t,\cdot)$ est inclus strictement dans une boule g\'eod\'esique de $M^m$. En utilisant l'unicit\'e de la g\'eod\'esique 
 reliant deux points d'une telle boule et la continuit\'e dans l'espace des chemins de ces g\'eod\'esiques par rapport \`a ces deux points pour construire une homotopie d\'eformant
 $u_\ep(s,\cdot)$ en une suite de chemins constants. Donc $[u_\ep(s,\cdot)]\in$ Im$\,ev_\ast$ mais comme $\pi_{k-1}(M^m)=0$ on obtient que $[u_\ep(t,\cdot)]=0$ ce qui contredit notre hypoth\`ese sur $u_\ep$.
 
 Afin de montrer que $W$ est effectivement la longueur d'une g\'eod\'esique ferm\'ee, on d\'eveloppe des outils qui permettent de d\'eformer les ensembles de niveau de E de fa\c con
 continue \`a un ``rythme donn\'e'' en l'absence de points critique. C'est la th\'eorie  des d\'eformations en dimension infine dite de {\bf Palais Smale}. Dans le cas pr\'esent son application est garantie par les faits suivants :
 \begin{itemize}
 \item{} L'{\bf espace de configuration} ${\mathcal M}:=W^{1,2}(S^1,M^m)$ d\'efinit une vari\'et\'e de Hilbert  r\'eguli\`ere, s\'eparable et model\'ee sur un espace de Hilbert $W^{1,2}(S^1,{\R}^m)$ lui m\^eme s\'eparable.
 \item{} ${\mathcal M}$ est \'equip\'e d'une m\'etrique
 \[
 g(u,v):=\int_{0}^{2\pi} u\cdot v+\p_\theta u\cdot\p_\theta v\ d\theta
 \]
 pour laquelle la variet\'e riemannienne associ\'ee $({\mathcal M},g)$ est \underbar{compl\`ete}.
 \item{} La fonctionelle $E$ est de r\'egularit\'e $C^1$ sur ${\mathcal M}$.
 \item{} La fonctionnelle $E$ satisfait la condition de {\bf Palais Smale}\footnote{Aussi appel\'ee ``condition (C)'' dans les textes d'origines.}  : Soit
 \[
 u_{n}\in {\mathcal M}\quad {t. q.}\quad\limsup_{n\rightarrow +\infty} E(u_n)<+\infty\quad \mbox{ et }\quad \nabla E_{u_n}\rightarrow 0
 \]
 alors il existe $u_{n'}$ et $u_\infty\in {\mathcal M}$ tels que
 \[
 u_{n'}\rightarrow u_\infty\quad\mbox{ et }\quad \nabla E_{u_\infty}=0\quad.
 \]
 \end{itemize}
Ces hypoth\`eses \'etant v\'erifi\'ees on d\'emontre sans trop de difficult\'ees  que le flot  $\phi_t$ associ\'e au champs de vecteur $X:=-\nabla E$ sur ${\mathcal M}$ existe \underbar{pour tout temps}. 
S'il n'y avait pas de point critique de $E$ entre les niveau d\'energie $W+\ep$ et $W-\ep$, la condition de Palais Smale nous dirait que $|\nabla E|$ est major\'ee par en dessous par un nombre strictement positif $\delta$ sur $E^{-1}(W-\ep,W+\ep)$. On prend alors un $u\in  C^0(S^{k-1},{\mathcal M})$ tel que $[u]\ne 0$ et satisfaisant
\[
\max_{s\in S^{k-1}}E(u(s,\cdot))<W+\ep\quad.
\]
 L'application du flot \`a toute la famille $(u(s,\cdot))_{s\in S^{k-1}}$ simultan\'ement nous fait passer sous le niveau d\'energie $W-\ep$ apr\`es le temps $T_0=2\ep/\delta$ :
 \[
 \max_{s\in S^{k-1}}E(\phi_{T_0}(u(s,\cdot)))<W-\ep\quad.
 \]
 Or, la continuit\'e du flot par rapport au param\`etre de temps donne $[\phi_{T_0}(u)]=[u]\ne 0$. Ce qui donne une condradiction et $W$ est bien la longueur d'une g\'eod\'esique ferm\'ee. On a ainsi 
 \'etabli le th\'eor\`eme suivant d\^u dans sa forme d\'efinitive la plus g\'en\'erale \`a Fet et Lyusternik.
 
 \begin{theo}
 \label{FL}\cite{FL}
 Toute vari\'et\'e riemannienne ferm\'ee poss\`ede une g\'eod\'esique ferm\'ee non constante.
 \end{theo}
 
 Dans son article d'origine datant de 1917 Birkhoff n'avait \'evidemment pas \`a sa disposition cette th\'eorie de la d\'eformation qui date de la deuxi\`eme moiti\'e des ann\'ees 60. N\'eanmoins
 il con\c cut ``\`a la main'', par un argument combinatoire, un processus de d\'eformation des ensembles de niveau d'\'energie, le {\bf processus de raccourcissement des courbes} qui peut \^etre vu comme un anc\^etre de ce que l'on appelle aujourd'hui les {\bf flots pseudo-gradients}, qui eux-m\^eme g\'en\'eralisent $X=-\nabla E$ au cas  de  fonctionnelles $C^1$ sur des vari\'et\'es non plus Hilbertiennes mais de Banach  munies de structures Finslerienne compl\`etes. L'utilisation des flots pseudo-gradients \`a partir du d\'ebut des ann\'ees 70  a elle-m\^eme stimul\'e l'\'etude des flots g\'eom\'etriques, qui en sont des versions plus pr\'ecises, et qui ne cessent de prouver leur efficacit\'e pour r\'esoudre diverses questions d'analyse, de g\'eom\'etrie et de topologie diff\'erentielle.
  
\section{L'infinit\'e de g\'eod\'esiques dans une vari\'et\'e ferm\'ee} 

La r\'esolution donn\'ee ci-dessus de l'existence d'une g\'eod\'esique sur une vari\'et\'e ferm\'ee quelconque pose naturellement la question du nombre de ces g\'eod\'esiques.
Au vue de cette r\'esolution il est naturel de proc\'eder \`a des op\'erations de minmax de la forme  (\ref{1}) en contraignant la classe $[u]$ \`a appartenir \`a des  groupes d'homotopie non triviaux \'eventuellement de plus en plus nombreux. La th\'eorie des d\'eformations de Palais-Smale nous permet d'affirmer que chaque  valeur de minmax est effectivement r\'ealis\'ee par une g\'eod\'esiques ferm\'ees. Se pose n\'eanmoins la question de la ``r\'ep\'etition'' de ces g\'eod\'esiques et il s'agit de savoir si la g\'eod\'esique obtenue n'est pas le recouvrement multiple d'une autre obtenue pr\'ec\'edemment. Une g\'eod\'esique ferm\'ee qui n'est pas le recouvrement d'une g\'eod\'esique ferm\'ee strictement plus courte est appel\'ee {\bf g\'eod\'esique premi\`ere}.

Dans le cas \`a nouveau o\`u le groupe fondamental de la vari\'et\'e $M^m$ poss\`ede une infinit\'e de classes de conjugaisons premi\`eres ({\it i.e.} qui ne peuvent \^etre repr\'esent\'ees par la composition it\'er\'ee d'une m\^eme classe) alors l'argument de minimisation mentionn\'e plus haut donne l'infinit\'e de g\'eod\'esiques primitives\footnote{Dans le cas de surfaces compactes hyperboliques $M^m=\Sigma_g$ de genre $g\ge 2$, comme la courbure est n\'egative, un r\'esultat classique  affirme que chaque classe du $\pi_1$ poss\`ede exactement une g\'eod\'esique. En effet, le calcul  de la d\'eriv\'ee seconde de $E$ donne imm\'ediatement  que chaque g\'eod\'esique est stable. Une \'etude approfondie des classes primitives de $\pi_1(\Sigma_g)$ permet alors d'\'etablir que le nombre de g\'eod\'esiques primitives de longueur major\'ee par $L>0$ se comporte asymptotiquement comme $e^L/L$. Ce r\'esultat de Delsarte, Huber et Selberg est aussi connu sous le nom du {\bf``th\'eor\`eme des nombres premiers des surfaces de Riemann''} du fait de la similitude avec la r\'epartition asymptotiques de ceux l\`a (voir par exemple \cite{Bus}). Il est encore bien plus d\'elicat de faire une estimation asymptotique des {\bf g\'eod\'esiques simples} - plong\'ees - parmi ces g\'eod\'esiques primitives. C'\'etait l'objet de la th\`ese de Maryam Mirzakhani qui a obtenu un comportement de la forme $\simeq n(\Sigma_g)\, L^{6g-6}$ \cite{Mir}.}.

Il est donc int\'er\'essant de se pencher plus sp\'ecifiquement sur le cas des vari\'et\'es dont le groupe fondamental est fini ou m\^eme nul. Pour le choix d'une {\bf m\'etrique g\'en\'erique}\footnote{ Un sous ensemble d'un espace topologique est ``g\'en\'erique'' s'il est l'intersection au plus d\'enombrable d'ensembles ouverts. Dans le cas d'espaces m\'etriques complets, par exemple, le th\'eor\`eme de Baire dit qu'un tel ensemble est dense. Une m\'etrique ``g\'en\'erique'' aussi appel\'ee ``cahoteuse'' - bumpy en anglais - sera suppos\'ee g\'en\'erique pour la topologie $C^3$. } th\'eor\`eme dit ``th\'eor\`eme des m\'etriques cahoteuses'' du \`a Ralph Abraham \cite{Abr} et Dmitry Anosov \cite{Ano}  affirme que pour un ensemble de m\'etriques g\'en\'eriques, toutes les g\'eod\'esiques de  l'espace des chemins ferm\'es ${\mathcal M}:=W^{1,2}(S^1,M^m)$, n'ont comme champs de Jacobi ({\it i.e.} \'elements du noyau de la d\'eriv\'ee seconde de $E$) que ceux g\'en\'er\'es par l'action du cercle $S^1$ sur le domaine. Comme les indices des points critiques de $E$ sont uniform\'ement born\'es sous un certain niveau d'energie et comme $E$ satisfait \`a la condition de Palais Smale, la th\'eorie de Morse nous enseigne que ${\mathcal M}:= W^{1,2}(S^1,M^m)$ est homotope \`a un complexe cellulaire limite\footnote{Voir l'appendice du livre de John Milnor \cite{Mil}.} de complexes cellulaires finis
donn\'es par les sous ensembles de niveaux
\[
{\mathcal M}^\la:=\lf\{u\in W^{1,2}(S^1,M^m)\ ;\ \sqrt{E(u)}\le \la\rg\}\quad.
\]
Par ailleurs, pour un choix de m\'etrique g\'en\'erique sur $M^m$ les points critiques non constants de $E$ r\'ealisent une sous vari\'et\'e unidimensionnelle sur laquelle $S^1$ agit et on peut aussi (grace \`a la forte compacit\'e de l'espace des g\'eod\'esiques de longueur born\'ees par un nombre fix\'e) garantir que les orbites appartiennent \`a des niveaux d'\'energie distincts les uns des autres.

 Dans l'espoir de mettre en \'evidence un nombre infini de g\'eod\'esiques primitives, plut\^ot donc que de consid\'erer les groups d'homotopie de l'espace des chemins ferm\'es ${\mathcal M}:=W^{1,2}(S^1,M^m)$, la th\'eorie de Morse nous invite \`a consid\'erer les diff\'erents groupes d'homologie associ\'es aux sous-ensemble de niveaux de $\sqrt{E}$ dans ${\mathcal M}$ :
\[
H_k({\mathcal M}^b,{\mathcal M}^a;{\Z})\quad \mbox{ o\`u }\quad a<b\ \mbox{  et }\ k\in {\N}\quad.
\]
Dans \cite{Gro1} Mikhael Gromov consid\`ere pour tout r\'eel positif  $\lambda>0$ le nombre entier d\'efinit par
\[
\mbox{dm}({\mathcal M}^\lambda):=\sup\lf\{l\in {\N}\ ;\ H_l({\mathcal M},{\mathcal M}^\lambda;{\Z})=0\ ,\ \forall\, l\le k \rg\}
\]
$\mbox{dm}({\mathcal M}^\lambda)$ ainsi d\'efini est aussi le nombre naturel maximal tel que toute application d'un polyh\`edre de dimension $l\le \mbox{dm}({\mathcal M}^\lambda)$ dans ${\mathcal M}$ se d\'eforme homotopiquement en une application dans
${\mathcal M}^\la$. On a donc en particulier
\be
\label{egalhom}
H_l({\mathcal M}^\lambda;{\Z})=H_l({\mathcal M};{\Z})\quad\mbox{ pour } l\le \mbox{dm}({\mathcal M}^\lambda)-1\quad.
\ee
On peut aussi d\'efinir le {\bf spectre non lin\'eaire}\footnote{Cette notion de {\bf spectre non-lin\'eaire}  est au centre de cet expos\'e et reviendra plus loin dans le cadre des surfaces minimales. Elle coincide avec la notion classique de spectre dans le cas lin\'eaire (voir \cite{Gro2}). En guise d'illustration, si on remplace ${\mathcal M}$ par les \'el\'ements de $W^{1,2}(M^m,{\R})$ de norme $L^2$ \'egale \`a 1 et  si on prend $E(u):=\int_{M^m}|du|^2 dvol_M$, le spectre $(\la_k)_{k\in {\N}}$ ainsi d\'efini est bien le spectre de l'op\'erateur de Laplace Beltrami sur $M^m$ o\`u chaque valeur du spectre est r\'ep\'et\'ee autant de fois que la dimension de l'espace propre correspondant. On peut donc tr\`es bien avoir $\lambda_k=\lambda_{k+1}$.}
$(\lambda_k)_{k\in {\N}}$ comme \'etant
\[
\lambda_k:=\sup\lf\{\lambda\in {\R}^+\quad;\quad \mbox{dm}({\mathcal M}^\lambda)<k\rg\}
\]
On a alors le th\'eor\`eme suivant.
\begin{theo}
\label{th-Gro1} \cite{Gro1}
Soit $(M^m,g)$ une vari\'et\'e  riemannienne compacte et de groupe fondamental fini, alors lorsque $k$ tend vers l'infini on a le comportement asymptotique suivant
\be
\label{2}
\frac{\log \la_k}{\log k}=1 +O\lf(\frac{1}{\log k}\rg)\quad.
\ee
\end{theo}
La preuve de la majoration $\mbox{dm}({\mathcal M}^\lambda)\le C \ \la$  est assez \'el\'ementaire dans le cas d'une m\'etrique g\'en\'erique et lorsque $M^m$ est simplement connexe. La d\'emonstration par minmax du th\'eor\`eme~\ref{FL} produit une g\'eod\'esique $\gamma$ qui ne peux \^etre un minimum local strict de la longueur. Comme la d\'eriv\'ee seconde n'a pour noyau que l'action de $S^1$ sur la param\'etrisation, la g\'eod\'esique $\gamma$ doit avoir un indice non nul.
On se convainc ais\'ement que le $n-$recouvrement  de cette g\'eod\'esique a un indice au moins \'egal \`a $n$. Par ailleurs,  le calcul direct de la d\'eriv\'ee seconde de ${E}$  montre que la  longueur contr\^ole lin\'eairement l'indice de Morse des points critiques. En consid\'erant les recouvrement entiers de $\gamma$ on obtient donc l'existence d'une famille de g\'eod\'esiques d'indice $k_n$ asymptotiquement comparable \`a  $n$. La th\'eorie de Morse (\cite{Mil}) nous enseigne\footnote{Pour une m\'etrique g\'en\'erique il y a au plus une g\'eod\'esique de longueur donn\'ee. } alors que cette g\'eod\'esique est responsable d'un changement d'homologie des sous-ensembles de niveau pour les longueurs $n\ L(\gamma)$ et que
$H_{k_n}({\mathcal M}^{n\,L(\gamma)+\ep},{\mathcal M}^{n\,L(\gamma)-\ep};{\Z})\ne 0$ pour $\ep$ suffisamment petit. Si $ \mbox{dm}({\mathcal M}^{n\,L(\gamma)+\ep})>k_n+1$ et si $ \mbox{dm}({\mathcal M}^{n\,L(\gamma)-\ep})>k_n$ , on doit alors avoir par d\'efinition
$H_{k_n}({\mathcal M},{\mathcal M}^{n\,L(\gamma)-\ep};{\Z})=0$ et $H_{k_n+1}({\mathcal M},{\mathcal M}^{n\,L(\gamma)+\ep};{\Z})=0$.
La suite de Mayer-Vietoris d'homologie relative
\[
H_{k_n+1}({\mathcal M},{\mathcal M}^{n\,L(\gamma)+\ep};{\Z})\longrightarrow H_{k_n}({\mathcal M}^{n\,L(\gamma)+\ep},{\mathcal M}^{n\,L(\gamma)-\ep};{\Z})\longrightarrow H_{k_n}({\mathcal M},{\mathcal M}^{n\,L(\gamma)-\ep};{\Z})
\]
apporte une contradiction \`a cette hypoth\`ese. Pour tout $n\ge 1$  on a donc $$\mbox{dm}({\mathcal M}^{n\,L(\gamma)})\le k_n+1\simeq n$$
ce qui donne\footnote{On utilise que  $\mbox{dm}({\mathcal M}^\lambda)$ est par d\'efinition une application croissante en $\lambda$.} $\mbox{dm}({\mathcal M}^\lambda)\le C \ \la$.

Une autre d\'emonstration  de la majoration avec pour seule hypoth\`ese que $\pi_1(M^m)$ est fini fait appel \`a un r\'esultat de Denis Sullivan \cite{Sul} sur le calcul du model minimal de $M^m$ donnant la cohomologie de de Rham de ${\mathcal M}$. Ce calcul \'etablit que, sous l'hypoth\`ese $\pi_1(M^m)$ fini, il existe une suite arithm\'etique $(k_n)_{n\in {\N}^\ast}=(p\,n+q)_{n\in {\N}^\ast}$ o\`u $p\in {\N}^\ast$ et $q\in{\N}$ telle que $H^{k_n}({\mathcal M};{\R})\ne 0$.

Dans le cas d'une m\'etrique g\'en\'erique \`a nouveau, le calcul direct de la d\'eriv\'ee seconde de ${E}$  donne que la  longueur contr\^ole lin\'eairement l'indice de Morse des points critiques. En d'autre termes, sous le niveau de longueur $\lambda$ il ne peut-y avoir que des g\'eod\'esique d'indice inf\'erieur \`a $C\,\la$ pour une constante $C>0$. {\bf L'in\'egalit\'e faible de Morse} (\cite{Mil}) donne alors que pour $C\,\la\le k$
\[
H^k({\mathcal M}^\la;{\R})=0\quad.
\]
Donc en particulier comme $H^{k_n}({\mathcal M};{\R})\ne 0$. (\ref{egalhom}) impose alors que pour $\la= C^{-1}\, k_n$ on ait $\mbox{dm}({\mathcal M}^\la)\le k_n$. Comme $k_n$ a une progression arithm\'etique on en d\'eduit  la majoration $\mbox{dm}({\mathcal M}^\lambda)\le C \ \la$.

La minoration est moins imm\'ediate et repose, dans le cas simplement connexe, d'une part sur la construction d'une d\'ecomposition cellulaire de ${\mathcal M}$ \`a partir d'une triangulation finie de $M^m$, et  d'autre part sur l'existence d'une application, homotope \`a l'identit\'e, lipschitzienne de $M^m$ dans lui m\^eme, qui contracte le squelette 1-dimensionnel de cette triangulation finie de $M^m$ en un point. Le cas o\`u $M^m$ \`a un groupe fondamental fini se ram\`ene
au cas simplement connexe en travaillant d'abord sur le rev\^etement universel\footnote{Voir aussi le chapitre 7 de \cite{Gr}.}.

La suite exacte d'homologie pour les paires $({\mathcal M}^a,{\mathcal M}^b)$  pour $\la\in {\N}$, $b>a> C\ \la$  o\`u $C> \lf(\liminf_{\la\rightarrow +\infty}{\mbox{dm}({\mathcal M}^\lambda)}/{\la}\rg)^{-1}$   donne pour $\lambda$ assez grand
\be
\label{3}
 \forall \ i\le \la\quad\quad H_i({\mathcal M}^b,{\mathcal M}^a;{\Z})=0\quad.
\ee
Comme, pour le choix g\'en\'erique de m\'etrique effectu\'e, il y a au plus une orbite $S^1$ de g\'eod\'esiques \`a un niveau d'\'energie donn\'ee, { l'in\'egalit\'e faible de Morse}  combin\'ee avec la minoration dans (\ref{2}) donne que pour avoir un indice de Morse $k$ pour une longueur de g\'eod\'esique $\la$ on a d'une part dm$({\mathcal M}^\la)\le k$, mais aussi dm$({\mathcal M}^\la)\ge C\ \la$. On en d\'eduit que
\[
\mbox{Indice de Morse}(\gamma)\ge C\, L(\gamma)\quad.
\]
En appliquant \`a nouveau l'in\'egalit\'e faible de Morse on obtient ainsi l'existence d'une constante $C>0$ pour laquelle
l'in\'egalit\'e suivante est vraie
\[
\mbox{Card}\lf\{\mbox{g\'eod. de longueur}\le \la\rg\}\ge\sum_{i=0}^{[C\ \la]}\mu_i\ge \  \sum_{i=0}^{[C\, \la]}dim(H_i({\mathcal M};{\R}))
\]
o\`u $\mu_i$ est le nombre d'orbite $S^1$ dans ${\mathcal M}$ de g\'eod\'esiques d'indice de Morse $i$. Toute g\'eod\'esique est le recouvrement multiple de g\'eod\'esiques premi\`ere.
On obtient donc la minoration
\be
\label{4}
\mbox{Card}\lf\{\mbox{g\'eod. {\bf premi\`eres} de longueur}\le \la\rg\}\ge\ C\ \frac{\ds \sum_{i=0}^{[C\, \la]}dim(H_i({\mathcal M};{\R}))}{\la}\quad.
\ee
Cette minoration a \'et\'e am\'elior\'ee par Werner Ballmann and Wolfgang Ziller qui \'etablissent dans \cite{BaZ}, \`a nouveau pour des m\'etriques g\'en\'eriques sur une vari\'et\'e ferm\'ee simplement connexe et $\lambda$ suffisamment grand,
\[
\mbox{Card}\lf\{\mbox{g\'eod. {\bf premi\`eres} de longueur}\le \la\rg\}\ge\ C\ \max_{i\le [C\, \la]}dim(H_i({\mathcal M};{\R}))\quad.
\]
Si  la suite des  nombres de Betti de l'espace des chemins ferm\'es n'est pas uniform\'ement born\'ee, on obtient ainsi l'existence d'une infinit\'e de g\'eod\'esiques premi\`eres.
On retrouve de cette fa\c con, dans le cas particulier de m\'etriques g\'en\'eriques, le fameux th\'eor\`eme de Detlef Gromoll et Wolfgang Meyer.
\begin{theo}
\label{th-Gro-Mey} \cite{GrMe}
Soit $(M^m,g)$ une vari\'et\'e riemannienne ferm\'ee dont le groupe fondamental est fini et telle que la suite des nombres de Betti de l'espace des chemins ferm\'es ${\mathcal M}:=W^{1,2}(S^1,M^m)$ ne soit pas uniform\'ement major\'ee, alors $(M^m,g)$ poss\`ede un nombre infini de g\'eod\'esiques ferm\'ees et premi\`eres.
\end{theo}

La condition de Gromoll-Meyer est par exemple satisfaite si $M^m$ est simplement connexe et si l'ag\`ebre  cohomologique $H^\ast(M,{\R})$ n'est pas g\'en\'er\'ee par un \'el\'ement unique, c'est en fait \'equivalent au fait que la suite $dim(H_i({\mathcal M};{\R}))$ ne soit pas born\'ee. Ce r\'esultat a \'et\'e d\'emontr\'e dans le cadre du calcul des mod\`eles minimaux de Sullivan pour les espaces de chemins ferm\'es dans \cite{PoSu}\footnote{Voir aussi une d\'emonstration dans \cite{FOT}.}

Pour conclure cette premi\`ere partie de l'expos\'e sur les g\'eod\'esiques ferm\'ees, notons que lorsque la vari\'et\'e $M^m$ est une vari\'et\'e simplement connexe quelconque, par exemple dans le cas o\`u $M^m$ est diff\'eomorphe \`a $S^m$, le probl\`eme de l'existence d'une infinit\'e de g\'eod\'esiques primitives est toujours ouvert si $m>2$. Le cas de la sph\`ere bidimensionnelle  a lui \'et\'e r\'esolu. L'on d\'emontre qu'il existe une infinit\'e de g\'eod\'esiques ferm\'ees premi\`eres sur $S^2$ en combinant les travaux de Victor Bangert \cite{Ban} et de John Franks \cite{Fra}.

\section{Existence de surfaces minimales ferm\'ees : l'approche param\'etrique}

Dans leur r\'esolution du probl\`eme de Plateau\footnote{Le probl\`eme de Plateau qui tient son nom du physicien Joseph Plateau  consiste \`a trouver un disque d'aire minimale bordant une courbe simple donn\'ee de ${\R}^3$.}, Jesse Douglas et Tibor Rad\'o minimisent l'\'energie de Dirichlet 
\[
E(u):=\frac{1}{2}\int_{D^2}|\nabla u|^2\ dx^2
\]
sur l'ensemble des applications dans l'espace de Sobolev $W^{1,2}(D^2,{\R}^3)$ dont la trace r\'ealise une fonction monotone continue du bord du disque $\p D^2$ dans la courbe de Jordan fix\'ee $\Gamma$ - cette condition de bord est  appel\'ee {\it condition de Plateau}. La raison pour laquelle la minimisation de l'\'energie de Dirichlet, plut\^ot que l'aire elle m\^eme qui est d'un usage variationnel tr\`es ``inconfortable'', permet de conclure vient essentiellement du fait que l'aire de $u$ est major\'ee par l'\'energie de $u$ avec \'egalit\'e si  et seulement si $u$ est faiblement conforme. Le theor\`eme d'uniformisation\footnote{Il s'agit plus pr\'ecisemment de la  version  approxim\'ee du th\'eor\`eme d'uniformisation pour les applications dans $W^{1,2}(D^2,{\R}^3)$ (\cite{Mor} chapitre I section 9).} permet de s'y ramener et un minimum de $E$ sous les contraintes de bord de Plateau sera automatiquement en param\'etrisation conforme et r\'ealisera un minimum de l'aire sous ces m\^eme contraintes.

Dans \cite{SaU}, Jonathan Sacks et Karen Uhlenbeck  reprennent, en dimension 2 cette fois,  la probl\'ematique d'Hadamard-Poincar\'e  sur l'existence de minimum de l'aire dans une classe du $\pi_2(M^m)$ donn\'ee.  La minimisation de l'\'energie de Dirichlet dans une classe d'homotopie de $\pi_2(M^m)$  d\'eg\'en\`ere en un ``bouquet'' de sph\`eres minimales qui reli\'ees convenablement entre elles par des tubes bien choisis permettent de retrouver la classe de d\'epart. En fait, avec cette approche, seule une famille de g\'en\'erateurs du groupe d'homotopie libre de $\pi_2(M^m)$ - apr\`es quotient par l'action de conjugaison du $\pi_1(M^m)$ sur $\pi_2(M^m)$ - est a priori r\'ealisable par des sph\`eres minimales. En particulier, gr\^ace au th\'eor\`eme d'Hurewicz, on en d\'eduit que si $\pi_1(M^m)=0$, alors un ensemble de g\'en\'erateurs de $H_2(M^m,{\Z})$ est r\'ealis\'e par des sph\`eres minimales. Il est important de noter \`a ce stade que ces surfaces minimales sont des immersions qui ont \'eventuellement des points de branchements isol\'es. L'article de Sacks et Uhlenbeck a donn\'e naissance \`a tout un domaine de l'analyse invariante conforme et a th\'eoris\'e les ph\'enom\`enes connus actuellement sous le nom de {\bf concentration-compacit\'e}.

Si on se penche maintenant sur la question de construire des sph\`eres minimales d'indice non nul,   les strat\'egies de minmax d\'evelopp\'ees plus haut dans le cas 1-dimensionnel   se heurtent en dimension 2 \`a deux obstructions fondamentales. 
\begin{itemize}
\item[- ] L'\'energie de Dirichlet en dimension 2 ne v\'erifie plus la condition de Palais Smale.
\item[ -] L'espace $W^{1,2}(S^2,M^m)$ ne constitue pas une vari\'et\'e de Banach r\'eguli\`ere finslerienne.
\end{itemize}
On ne peut donc plus a priori appliquer directement la th\'eorie de d\'eformation en dimension infinie de Palais et Smale. Une alternative consiste alors \`a travailler dans l'espace plus petit $W^{1,2+\alpha}(S^2,M^m)$ qui lui d\'efinit bien une vari\'et\'e de Banach r\'eguli\`ere et \'equip\'ee d'une structure finslerienne compl\`ete. L'\'energie ``renforc\'ee''\footnote{Cette \'energie renforc\'ee $E_\al$ a en fait \'et\'e introduite \`a l'origine pour le probl\`eme de minimisation dans les classes d'homotopie par Sacks et Uhlenbeck.} 
\[
E_\al(u):=\frac{1}{2}\int_{S^2}\lf(1+|d u|_{S^2}^2\rg)^{1+\al}\ dvol_{S^2}
\]
v\'erifie bien la condition de Palais Smale et, apr\`es avoir exhib\'e des points critiques de $E_\al$ pour une valeur de minmax donn\'ee, on peut faire tendre $\al$ vers 0 et l'analyse de Sacks et Uhlenbeck permet d'obtenir
un ``bouquet'' de points critiques de l'\'energie de Dirichlet mais dont il faut s'assurer que les applications associ\'ees sont conformes afin de pouvoir affirmer que l'ensemble correspond bien \`a un bouquet de sph\`eres minimales.

La question de l'existence de surfaces minimales par l'approche param\'etrique se complique si on cherche \`a prescrire le genre de cette surface. Dans le cas $S^2$, comme il n'existe qu'une seule classe conforme, la minimisation de l'\'energie de Dirichlet garantit qu'un minimum est bien conforme et donc r\'ealise une sph\`ere minimale. Ce n'est plus le cas si on remplace $S^2$ par une autre surface de Riemann ferm\'ee $\Sigma$ de genre non nul. Lors d'une proc\'edure de minimisation de $E$ il faut aussi introduire la possibilit\'e pour la classe conforme sur $\Sigma$ de parcourir ``librement'' l'espace de Teichm\"uller associ\'e. L'espace de configuration n'est plus alors $W^{1,2}(\Sigma, M^m)$ mais le produit
$W^{1,2}(\Sigma,M^m)\times {\mathcal T}_\Sigma$ o\`u ${\mathcal T}_\Sigma$ d\'esigne l'espace de Teichm\"uller de $\Sigma$. Le contr\^ole de la classe conforme est \'evidemment une difficult\'e suppl\'ementaire de taille lors d'un processus de minimisation. Dans certaines situations n\'eanmoins, cette difficult\'e peut \^etre lev\'ee. C'est le cas par exemple si on se restreint \`a des applications $u$ qui envoient injectivement
le $\pi_1$ de $\Sigma$ sur le $\pi_1$ de $M^m$. De telles applications sont dites ``incompressibles''. Sous l'hypoth\`ese d'incompressibilit\'e on d\'emontre alors que l'\'energie de Dirichlet 
\[
E_h(u):=\frac{1}{2}\int_\Sigma|du|^2_h \, dvol_h
\]
contr\^ole la classe conforme $h$ de $\Sigma$, repr\'esent\'ee ici par une m\'etrique de courbure scalaire constante de volume 1. Une cons\'equence de la pr\'ecompacit\'e de la classe conforme dans l'espace des modules sous l'hypoth\`ese du contr\^ole de l'\'energie  est l'existence d'une application $u$ harmonique de $(\Sigma,h)$ dans $M^m$ pour une surface de Riemann $(\Sigma,h)$ dans cet espace d'applications incompr\'essibles qui est conforme et qui minimise l'aire dans cet espace.

Chaque surface incompressible g\'en\`ere ce que l'on appelle
un {\bf sous groupe de surface} qui est l'image (injective) par $u_\ast$ de $\pi_1(\Sigma)$  dans  $\pi_1(M^m)$ (cette image est par d\'efinition isomorphe au $\pi_1(\Sigma)$).

Dans \cite{KaM}, Jeremy Kahn et Vladimir Markovi\'c d\'emontrent que  pour  une vari\'et\'e hyperbolique ferm\'ee donn\'ee $(M^3,g)$, il existe au moins $(c_1\, \gamma)^{2\gamma}$ classes de conjugaisons de sous groupes de surfaces de genre $\gamma$. Les surfaces $\Sigma$, tout comme la vari\'et\'e hyperbolique choisie $(M^3,g)$, sont des espaces asph\'eriques\footnote{Les espaces asph\'eriques sont les espaces dont les rev\^etements universels  sont contractiles et  leurs seuls groupes d'homotopies non triviaux sont leurs $\pi_1$.}. Pour cette raison, deux applications d\'efinissant deux classes de conjugaison distinctes pour l'homomorphisme $u_\ast$ : $\pi_1(\Sigma)\rightarrow \pi_1(M^3)$ sont aussi dans des classes distinctes de l'espace des classes d'homotopies libres $[\Sigma, M^3]$ entre $\Sigma$ et $M^3$. En d'autres termes, chaque classe de conjugaison  des sous groupes de surfaces de genre $\gamma$ donn\'e correspond \`a des composantes connexes distinctes de l'espace des applications incompr\'essibles. Ainsi, en utilisant l'analyse d\'ecrite ci dessus,  on peut associer \`a chacune de ces classes  une surface minimale de genre $\gamma$. Tout comme dans le cas des g\'eod\'esiques, afin de ne prendre en compte que les recouvrements simples de surfaces minimales, il faut dans le cas pr\'esent garder \`a l'esprit qu'une surface de Riemann peut \^etre le recouvrement non branch\'e d'une autre surface de genre strictement sup\'erieur. L'application de recouvrement par ailleurs induit naturellement une injection du groupe fondamental. Un sous groupe de surface dans une vari\'et\'e hyperbolique ferm\'ee $(M^3,g)$ est dit {\bf premier} s'il n'est pas r\'ealis\'e par la composition successive d'un recouvrement non branch\'e sur une surface de genre inf\'erieur suivi de l'image de cette surface par une application incompressible. Comme la croissance du nombre de groupe de surfaces est au moins en $(c_1\, \gamma)^{2\gamma}$ la grande majori\'e 
de ces groupes de surfaces sont premiers dans $M^3$ et on a d\'emontr\'e le r\'esultat suivant.
\begin{theo}
\label{hyperb}\cite{KaM}
Sur une vari\'et\'e  hyperbolique ferm\'ee donn\'ee il existe une infinit\'e d'immersion\footnote{Comme ces immersion minimisent l'aire dans la classe d'homotopie libre donn\'ee, un r\'esultat de Osserman \cite{Oss} et Gulliver \cite{Gul}  afirme que ces surfaces minimales en dimension 3 sont exemptes de points de branchements.} minimales g\'eom\'etriquement distinctes.
\end{theo}
Il est int\'erressant d'observer \`a ce stade que la famille infinie de surfaces produites par Kahn et Markovi\'c est faite de surfaces qui sont toutes stables (en fait qui minimisent l'aire). Ceci contraste avec la famille infinie que nous allons construire dans la derni\`ere partie de cet expos\'e o\`u les surfaces seront d'indice de Morse tendant vers l'infini.
\section{Les variations de l'aire en th\'eorie de la mesure g\'eom\'etrique}\footnote{Le lecteur non familier avec la th\'eorie de la mesure g\'eom\'etrique est invit\'e \`a se reporter \`a l'expos\'e pr\'ec\'edent \cite{Riv}.}

Le but de la section pr\'ec\'edente \'etait principalement d'illustrer la difficult\'e que constitue l'extension en dimension sup\'erieure \`a  1 de l'approche param\'etrique 
pour les op\'erations de minmax, et plus pr\'ecisemment de l'utilisation pour les surfaces minimales de la th\'eorie des d\'eformations de Palais et Smale.
En fait, le constat de la limitation de l'approche param\'etrique \`a des fins variationnelles pour les surfaces minimales \'etait d\'ej\`a pr\'esent bien ant\'erieurement aux travaux
d'analyse globale de la fin des ann\'ees 60. La th\'eorie de la mesure g\'eom\'etrique est n\'ee - au milieu des ann\'ees 50 - du projet d'\'etendre le probl\`eme de Plateau aux dimensions sup\'erieures avec
 des donn\'ees plus  g\'en\'erales que des courbes de Jordan dans ${\R}^3$. Dans sa premi\`ere phase l'objet principal de la TMG \'etait la notion de {\bf courant r\'ectifiable entiers}, fusion entre les courants de de Rham (distributions vectorielles) et le concept de rectibiabilit\'e\footnote{Pour $0\le k\le Q$, un sous ensemble de ${\R}^Q$ est dit {\bf d\'enombrablement $k-$rectifiable} si, en dehors d'un ensemble de mesure de Hausdorff $k-$dimensionel nulle, il  est contenu dans une union au plus d\'enombrable d'images par des fonctions lipschitzienne de sous ensembles de ${\R}^k$ et sa ${\mathcal H}^k-$masse est finie.  } de Besicovitch. Le th\'eor\`eme fondateur de cette th\'eorie est celui d'Herbert Federer et de Wendell Fleming qui \'etablit que l'espace des courants rectifiables entiers de bord nul - ou de bord de masse\footnote{La masse d'un courant est une quantit\'ee non n\'egative qui coincide avec le volume lorsque le courant est par exemple donn\'e par une sous-vari\'et\'e.} contr\^ol\'ee - dans une vari\'et\'e riemannienne quelconque et de masse uniform\'ement born\'ee est faiblement s\'equentiellement ferm\'e\footnote{En fait cette compacit\'e faible s\'equentielle est valable pour la topologie dite ``b\'emol'',  une topologie un peu plus fine que celle donn\'ee par la simple dualit\'e avec les formes diff\'erentielles r\'eguli\`ere \`a support compact que nous introduisons plus loin dans l'expos\'e. L'auteur de ces notes invite le lecteur  \`a se reporter par exemple au s\'eminaire Bourbaki no 1080 \cite{Riv} pour une description un peu plus d\'etaill\'ee de certains fondamentaux de la th\'eorie de la mesure g\'eom\'etrique.}. Comme la masse est faiblement sequentiellement sous-continue, un corollaire de ce r\'esultat est l'existence d'un courant rectifiable entier de masse minimale dans toute classe d'homologie enti\`ere fix\'ee (relative ou pas relative). Un tel courant est appel\'e {\bf courant rectifiable entier minimisant l'aire}. 
 
 Un tel objet a-priori pourrait \^etre  tr\`es singulier (c'est le prix \`a payer pour s'\^etre plac\'e dans un espace faiblement sequentiellement ferm\'e pour un contr\^ole de la masse - La masse en g\'eom\'etrie est une quantit\'e tr\`es peu coercive). Il n'en est rien, en particulier dans le cas de la codimension 1 qui va nous occuper exclusivement d\'esormais jusqu'\`a la fin de cet expos\'e. Le th\'eor\`eme suivant est le r\'esultat d'efforts combin\'es de De Giorgi, Triscari, Fleming, Simons, Federer, Almgren.
 
 \begin{theo}
 \label{reg}\cite{DG}, \cite{Fle}, \cite{Sims}, \cite{Fed}, \cite{Almg}
  Soit $(M^m,g)$ une vari\'et\'e riemannienne et soit  $T$ un courant rectifiable entier sans bord minimisant l'aire dans une classe d'homologie relative  $H_{m-1}(\Om,\p\Om;{\Z})$ o\`u $\Om$ est un ouvert quelconque de $M^m$. Alors $T$ est le courant d'int\'egration d'une sous-vari\'et\'e plong\'ee de $\Om$ en dehors d'un sous ensemble ferm\'e de codimension de Hausdorff 8 et \'equip\'e d'une multiplicit\'e enti\`ere r\'eguli\`ere. En dimension 
 8 l'ensemble singulier est inclus dans un ensemble de points isol\'es.
 \end{theo}
Ce r\'esultat est optimal : Le c\^one issue de la sous-vari\'et\'e de la sph\`ere de dimension 7 donn\'ee par $S^3_{1/\sqrt{2}}\times S^3_{1/\sqrt{2}}$ minimise l'aire \cite{BGG} (fig. 1). Cet \'etat de fait est fondamental dans notre expos\'e et se trouve \^etre derri\`ere la limitation en dimension du th\'eor\`eme principal qui nous int\'eresse aujourd'hui.

\begin{center}
\begin{tabular}{lll}
\footnotesize{\textcolor{blue}{$\partial C_3 =$} \textcolor{red}{$S^1_{\frac{1}{\sqrt{2}}} \times S^1_{\frac{1}{\sqrt{2}}}$} $\subset S^3 \subset \mathbb{R}^4$~~~~~}
&
\footnotesize{\textcolor{blue}{$\partial C_5 =$} \textcolor{red}{$S^2_{\frac{1}{\sqrt{2}}} \times S^2_{\frac{1}{\sqrt{2}}}$} $\subset S^5 \subset \mathbb{R}^6$}
&
\footnotesize{\textcolor{blue}{~~~~~$\partial C_7 =$} \textcolor{red}{$S^3_{\frac{1}{\sqrt{2}}} \times S^3_{\frac{1}{\sqrt{2}}}$} $\subset S^7 \subset {\R}^8$}
 \end{tabular}
\end{center}
\begin{center}
\psfrag{C3}{\textcolor{blue}{$C_3$}}
\psfrag{C5}{\textcolor{blue}{$C_5$}}
\psfrag{C7}{\textcolor{blue}{$C_7$}}
\psfrag{4D}{$4D$}
\psfrag{6D}{$6D$}
\psfrag{8D}{$8D$}
\psfrag{p}{\textcolor{green}{\footnotesize point}}
\psfrag{s}{\textcolor{green}{\footnotesize singulier}}
\psfrag{CS}{C\^one de Simons}
\includegraphics[width=14cm]{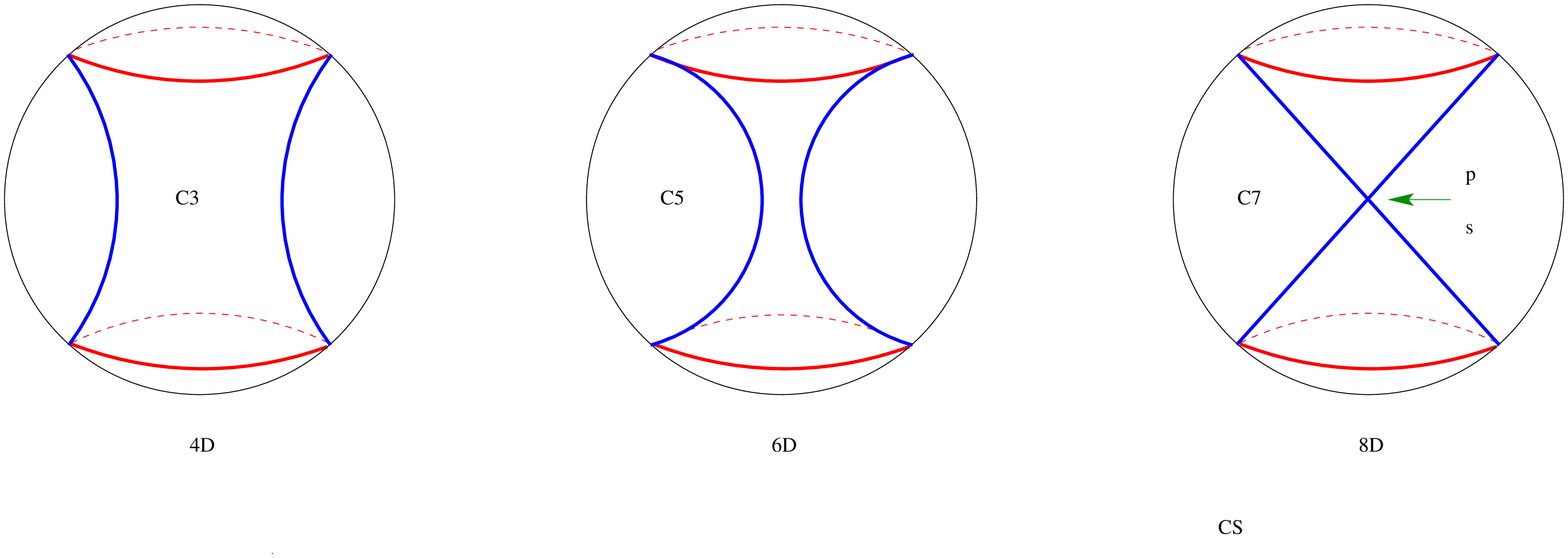}
\end{center}

\bigskip
\begin{center}
Fig.~1: ~ Hypersurfaces minimisant le volume
\end{center}

Si le cadre des courant entiers rectifiables est efficace pour mettre en \'evidence des hyper-surfaces minimales minimisant l'aire dans des classes d'homologie enti\`eres donn\'ees, ce n'est plus le cas d\`es lors que l'on cherche
\`a construire variationnellement des hyper-surfaces minimales d'indice non nul. La raison en est que la masse n'est que semi-continue inf\'erieure pour la convergence faible des courants et que dans une op\'eration
de type minmax il faut pouvoir avoir la continuit\'e de celle ci. Cette constatation a stimul\'e la ``seconde g\'en\'eration'' de la th\'eorie de la mesure g\'eom\'etrique sous l'impulsion des contributions succ\'essives de Fred Almgren \cite{Alm1}, William Allard \cite{All} et Jon Pitts \cite{Pitt}. L'objet principal  de cette deuxi\`eme partie de la th\'eorie devient une version - tr\`es faible ! - de type ``th\'eorie de la mesure'' d'une sous-vari\'et\'e il s'agit de la notion de {\bf varifold} introduite par Almgren\footnote{{\it ``I called the new objects ``varifolds'' having in mind that they were a measure theoretic substitute for manifolds created for the variational calculus'' \cite{Alm2}.}}. Un varifold de dimension $k$ dans une vari\'et\'e riemannienne $M^m$ n'est autre qu'une {\bf mesure de Radon dans  la vari\'et\'e de Grassmann des $k-$plans  non orient\'es du fibr\'e tangent \`a $M^m$}, $G_k(TM^m)$.  On peut par exemple adjoindre \`a toute sous-vari\'et\'e $k-$dimensionnelle $\Sigma$ de la vari\'et\'e riemannienne $M^m$, le {\it varifold} associ\'e qui est explicitement  donn\'e par 
\[
\forall\ \Phi\in C^0(G_k(TM^m))\quad\quad {\mathbf v}_{\Sigma}(\Phi):=\int_\Sigma \Phi(x,T_x\Sigma)\ dvol_\Sigma
\]
L'id\'ee derri\`ere la notion des varifolds est inspir\'ee de la notion, commune dans les probl\`emes de convergence faible en EDP, de {\bf mesures de Young} et consiste \`a ``enregistrer'' l'information donn\'ee par le premier jet lors de passages \`a la limite. La masse par ailleurs passe \`a la limite pour la convergence faible des mesures de Radon et c'est principalement pour cela que ce type d'objets a \'et\'e introduit. Les varifolds  de dimension $k$ sont appel\'es {\bf rectifiables} si leur projection sur l'espace des mesures de Radon de $M^m$ est un ensemble rectifiable de $M^m$ et pour ${\mathcal H}^k-$presque tout point de cet ensemble rectifiable la mesure sur la fibre est donn\'ee par la masse de Dirac du $k-$plan approxim\'e tangent. Un varifold rectifiable est dit {\bf stationnaire} s'il est point critique de la masse (la variation totale de la mesure de Radon correspondante) pour les perturbations infinit\'esimales de l'identit\'e dans l'espace des diffeomorphismes de $M^m$.

De la m\^eme fa\c con que dans les sections pr\'ec\'edents nous allons d\'efinir un {\bf espace de configuration} sur lesquels nous allons poser des probl\`emes de minmax sur l'aire (ou plut\^ot la ``masse'' en TMG) et effectuer des op\'erations de d\'eformation afin de d\'emontrer que les valeurs de minmax sont atteintes par des surfaces minimales.

Nous avons d\'ej\`a vu des espaces de configurations donn\'es par des espaces de Sobolev pour les constructions par minmax de g\'eod\'esiques o\`u de surfaces minimales dans les espaces hyperboliques tridimensionnels. L'espace que nous introduisons maintenant est un espace d'objets bien moins r\'eguliers, l'espace des {\bf chaines b\'emol \`a coefficient dans ${\Z}_2$}. Soit $G$ un groupe ab\'elien \'equip\'e d'une norme qui en fait un espace m\'etrique complet (on va consid\'erer esentiellement $G={\Z}$, ${\R}$ ou ${\Z}_p$). Soit $P_k({\R}^Q;G)$ 
l'espace des chaines polyh\'edrale\footnote{Une $G-$cha\^\i ne polyh\'erale est une somme finie  de polyh\`edres convexes born\'es de dimension $k$ de ${\R}^Q$, ne s'intersectants pas mutuellements et  affect\'e chacun d'une ``multiplicit\'e'' dans $G$. \`A une telle chaine $A:=\sum_{i=1}^ng_i\sigma_i$ on associe tout d'abord la masse ${\mathbf M} (A)=\sum_{i=1}^n|g_i|\ {\mathcal H}^k(\sigma_i)$ puis une application ${\mathcal H}^k$ mesurable, born\'ee \`a valeur dans $G$ \'egale \`a $g_i$ sur chaque $\sigma_i$ et \'egale \`a l'\'el\'ement neutre ailleurs. Finalement $P_k({\R}^Q,G)$ est quotient\'e par la relation d'\'equivalence suivante : deux \'el\'ements sont \'egaux s'ils d\'efinissent la m\^eme application ${\mathcal H}^k$ presque partout. On v\'erifie alors que l'op\'erateur de bord habituel sur les chaines est bien d\'efini c'est \`a dire qu'il est ind\'ependant du repr\'esentant choisi - voir \cite{Whit}.} \`a coefficient dans $G$.  Dans \cite{Whit} chapitre 5, Hassler Whitney introduit sur $P_k({\R}^Q;G)$ la norme b\'emol\footnote{Whitney introduit aussi une {\bf norme di\`ese} sur les chaines polyh\`edrales - born\'ee par la norme b\'emol et donc paradoxallement ''moins \'elev\'ee'' pour reprendre l'image musicale - \`a moins qu'il ne s'agisse des cochaines o\`u l\`a l'ordre des normes s'inverse... La compl\'etion  de  $P_k({\R}^Q;G)$ pour la norme di\`ese 
ne permet pas d'\'etendre continuement l'op\'erateur de bord contrairement \`a la norme b\'emol. Ce qui peux expliquer pourquoi la premi\`ere est moins connue que la seconde.} :
\[
{\mathcal F}(A):=\inf\lf\{M(A-\p D)+M(D)\ ; \ D\in {\mathcal P}_{k+1}({\R}^Q,G)\rg\}
\]
L'espace $I_k({\R}^Q;G)$ des {\bf $k-$chaines b\'emol \`a valeur dans $G$} est la completion de $P_k({\R}^Q;G)$ dans l'espace des $G-$courants . Un th\'eor\`eme important du \`a Fleming \cite{Fle1} (lorsque $G$ est fini) et White \cite{Whi1} (cas g\'en\'eral) affirme que tout \'el\'ement de masse finie dans $I_k({\R}^Q;G)$ est {rectifiable}. On note $\ti{I}_k({M}^m;G)$ l'ensemble des chaines b\'emol inclues dans $M^m$ de masse finie et dont le bord est de masse finie.

Les espaces de configuration qui vont nous inter\'esser pour la suite de cet expos\'e sont les {\bf hypercycles \`a coefficient ${\Z}_2$ } donn\'es par
\[
{\mathcal M}:=Z_{m-1}(M^m;{\Z}_2):=\lf\{ A\in \ti{I}_{m-1}(M^m;{\Z}_2)\ ;\ \partial A=0\rg\} 
\]
ou encore, lorsque $M^m$ n'est pas ferm\'e et est compacte avec bord on d\'efinit dans un premier temps
\[
Z_{m-1}(M^m,\p M^m;{\Z}_2):=\lf\{ A\in \ti{I}_{m-1}(M^m;{\Z}_2)\ ;\ \partial A\subset \p M^m\rg\} \quad.
\]
et on d\'efinit l'espace des {\bf hypercycles relatifs ${\Z}_2$} not\'e ${\mathcal M}:=Z_{m-1,rel}(M^m,\p M^m;{\Z}_2)$ comme \'etant le quotient de $Z_{m-1}(M^m,\p M^m;{\Z}_2)$ par la relation
\[
A\equiv B\quad\mbox{ si } \quad A-B\in \ti{I}_{m-1}(\p M^{m};{\Z}_2)
\]
Les efforts\footnote{Ce mot ``efforts'' cache en fait des travaux d'une difficult\'e  conceptuelle et technique tr\`es cons\'equente, le contexte ``mou'' des hypercycles compromet toute possibilit\'e d'impl\'ementation de la th\'eorie des d\'eformation de Palais-Smale. Nous avions essay\'e d'en rendre compte dans un expos\'e pr\'ec\'edent \cite{Riv}. Les applications continues \`a valeur dans les hypercycles ${\Z}_2$ doivent \^etre d\'eform\'es dans un premier temps en des applications continues pour des topologies plus fortes issue de la topologie des varifolds et de la masse. Elles sont dans un deuxi\`eme temps discr\'etis\'ees, puis, dans une proc\'edure inspir\'ee du raccourcissement des courbes de Birkhoff, ces applications sont ``s\'err\'ees'' - pull tight procedure - afin de mettre en \'evidence un hypercycle r\'ealisant la largeur et constituant un varifold stationnaire {\bf presque minimisant} sur les anneaux (la d\'efinition est due \`a Pitts \cite{Pitt}). Cette derni\`ere propri\'et\'e, qui s'apparente au fait d'\^etre localement minimisant et donc d'indice global fini, permit \`a Jon Pitts de d\'emontrer
la r\'egularit\'e de cet hypercycle en utilisant les travaux contemporains de Rick Schoen et Leon Simon \cite{ScS} \'etablissant des estim\'ees de courbure sur les surfaces stables minimales en codimension 1.}  combin\'es de Almgren  puis de Pitts permettent d'\'etablir le th\'eor\`eme suivant.
\begin{theo}
\label{almpitt} \cite{Alm1}, \cite{Pitt} Soit $M^m$ une vari\'et\'e ferm\'ee et $\Pi$ une classe non nulle de $\pi_k\lf(Z_{m-1}(M^m;{\Z}_2)\rg)$ alors la largeur correspondante
est strictement positive :
\[
W(\Pi):=\inf_{{\Phi}\in \Pi}\sup_{x\in S^k} {\mathbf M}(\Phi(x))>0
\]
et il existe un varifold stationaire $V$ port\'e par une sous-vari\'et\'e plong\'ee et r\'eguli\`ere en dehors d'un sous-ensemble ferm\'e de codimension 8 dans $M^m$  tel que ${\mathbf M}(V)=W(\Pi)$. En particulier,
si la dimension de $M^m$ est inf\'erieure o\`u \'egale \`a 7, il existe un nombre fini d'hypersurfaces minimales ferm\'ees et connexes, ne s'intersectant pas, $\Sigma_1\cdots \Sigma_n$ ainsi que $n$ entiers naturels non nuls
$q_1\cdots q_n$ tels que 
\[
V=q_1\, \Sigma_1+\cdots+q_n\, \Sigma_n\quad.
\]
\end{theo}

\section{Spectres non lin\'eaires et largeurs de Gromov.}
Afin de construire une surface minimale il suffit alors d'identifier au moins une classe d'homotopie non nulle dans ${\mathcal M}$. L'\'etude des groupes
d'homotopie des cycles b\'emol \'etait en fait l'objet de la th\`ese d'Almgren. Il d\'emontre\footnote{La preuve d'Almgren est \'ecrite dans le cas des cycles entiers et des homologies enti\`eres de $M^m$. Elle s\'etend cependant au cas de $G={\Z}_p$  - voir \cite{Pitt}.} en particulier le r\'esultat suivant.
\begin{theo}
\label{almhomot}\cite{Alm}
Soint $M^m$ une vari\'et\'e compacte connexe de dimension $m$. Pour tout $n\in {\N}$ et $k\in\{1\cdots m\}$ il existe un isomorphisme entre $\pi_n(Z_{k,rel}(M^m,\p M^m;{\Z}_2))$ et $H_{n+k}(M^m,\p M^m;{\Z}_2)$. En particulier la composante\footnote{Il s'agit en fait de la composante donn\'ee par les hypercycles qui sont des bords.} 
$Z^0_{m-1,rel}(M^m,\p M^m;{\Z}_2)$ de $Z_{m-1,rel}(M^m,\p M^m;{\Z}_2)$ contenant 0 poss\`ede exactement un groupe d'homotopie non nul $\pi_1(Z_{m-1,rel}(M^m, \p M^m;{\Z}_2))={\Z}_2$. 
\end{theo}
La deuxi\`eme partie de ce th\'eor\`eme est l'affirmation que $Z^0_{m-1,rel}(M^m,\p M^m;{\Z}_2)$ est un espace d'Eilenberg MacLane $K({\Z}_2,1)$ et donc qu'il existe une {\bf \'equivalence d'homotopie faible}
\[
Z^0_{m-1,rel}(M^m,\p M^m;{\Z}_2)\simeq_{faible} {\R}P^\infty
\]
o\`u ${\R}P^\infty$ est l'espace projectif infini (en effet la sph\`ere $S^\infty$ qui est contractible est le relev\'e universel et recouvre exactement 2 fois ${\R}P^\infty$).

Un g\'en\'erateur du $\pi_1(Z^0_{m-1,rel}(M^m,\p M^m;{\Z}_2))$ est donn\'e par un {\bf balayage de Birkhoff} de $M^m$ c'est \`a dire un chemin d'hypercycles de $M^m$ qui ``recouvre'' $M^m$ exactement une fois.

Cette \'equivalence faible d'homotopie est r\'ealis\'ee par l'application suivante. Soit $f$ une fonction de Morse sur $M^m$ \`a valeur dans $[0,1]$ (il en existe toujours \cite{Mil}). La famille donn\'ee par
\[
t\in [0,1]\ \longrightarrow\ \partial\lf(\lf\{ x\in M^m\ ;\ f(x)<t\rg\}\rg)\in Z^0_{m-1}(M^m, \p M^m;{\Z}_2)
\]
r\'ealise un balayage de Birkhoff et donc g\'en\`ere en particulier le $\pi_1$ des ${\Z}_2-$hypercycles. On l'\'etend \`a ${\R}P^\infty$ de la fa\c con suivante.
\[
\Phi([a_0,a_1\cdots, a_k,0,0\cdots]):=t\in [0,1]\ \longrightarrow\ \partial\lf(\lf\{ x\in M^m\ ;\ \sum_{i=1}^ka_i\  f^i(x)<t\rg\}\rg)\quad.
\]

En combinant le th\'eor\`eme~\ref{almpitt} et le th\'eor\`eme~\ref{almhomot} on obtient le r\'esultat suivant dont la d\'emonstration compl\`ete se trouve dans le livre de Pitts\footnote{Ce r\'esulat a \'et\'e prouv\'e plus ou moins en m\^eme temps en dimension 3 par Francis Smith dans sa th\`ese \cite{Smi} - non publi\'ee - \'ecrite sous la direction de Leon Simon. Le lecteur peut aussi se reporter \`a \cite{CoD} pour une d\'emonstration dans le cas tri-dimensionnel.}

\begin{coro}\cite{Pitt}
Soit $M^m$ une vari\'et\'e riemannienne ferm\'ee de dimension inf\'erieure ou \'egale \`a $7$, alors il existe dans $M^m$ une sous-vari\'et\'e minimale de codimension 1.
\end{coro}

Si maintenant on ambitionne d'\'etudier la possibilit\'e d'avoir plusieurs hypersurfaces minimales et \'eventuellement un nombre infini, les minmax construits autour des groupes d'homotopie ne suffisent pas\footnote{Il n'existe qu'un seul groupe d'homotopie non trivial.} . 
En analogie avec le cas des g\'eod\'esiques dans la premi\`ere partie de l'expos\'e, on pr\'ef\`ere alors une approche de type ``Morse'' consistant \`a \'etudier
les homologies des sous-ensembles de niveau de l'espace de configuration ${\mathcal M}= Z^0_{m-1,rel}(M^m, \p M^m,{\Z}_2)$. L'equivalence faible homotopique entre $Z^0_{m-1,rel}(M^m,\p M^m;{\Z}_2)$ et  ${\R}P^\infty$ induit un isomorphisme d'anneau cohomologiques\footnote{\`A ce stade on va pr\'ef\'erer travailler avec la cohomologie plut\^ot qu'avec l'homologie de $Z^0_{m-1}(M^m,{\Z}_2)$ car cette premi\`ere poss\`ede une structure d'anneau issue du cup-produit. L'espace des cha\^\i nes est probablement g\'eom\'etriquement plus intuitif et  de ce fait a priori plus agr\'eable \`a manipuler, mais la structure d'anneau de $H^\ast(\mathcal M;{\Z}_2)$ l'emporte.}  et on obtient finalement que l'anneau cohomologique de l'espace de configuration est un $\Z_2-$anneau polynomial, libre \`a un g\'en\'erateur de degr\'e un\footnote{On ne peux pas imaginer plus simple. Un espace de configuartion sur des objets plus lisses, et donc plus ``contraints'', que hypercycles ${\Z}_2$ aurait certainement induit un type d'homotopie plus complexe et plus riche mais aussi moins facile \`a calculer et \`a identifier au moyen d'op\'erations g\'eom\'etriques simples. C'est le grand avantage du point de vue que nous prenons \`a ce stade.}
\[
H^\ast(Z^0_{m-1,rel}(M^m,\p M^m;{\Z}_2))=\Z_2[\al]
\]
Le g\'en\'erateur $\al$ ``d\'etecte'' les balayages de Birkhoff au sens o\`u pour tout $$\Psi\in C^0_{\mathcal F}( S^1,Z^0_{m-1,rel}(M^m,\p M^m;{\Z}_2))$$  qui ``recouvre'' exactement une fois $M^m$ (et donc g\'en\`ere $\pi_1(Z^0_{m-1,rel}(M^m,\p M^m;{\Z}_2))$) on a
\[
\lf<\al,[\Psi]\rg>_{H^1({\mathcal M};{\Z}_2),H_1({\mathcal M},\p M^m;{\Z}_2)}\ne 0\quad.
\]
De fa\c con analogue au cas des g\'eod\'esiques dans la premi\`ere partie de l'expos\'e on introduit le {\bf spectre non lin\'eaire} associ\'e \`a la masse dans ${\mathcal M}:=Z^0_{m-1,rel}(M^m,\p M^m,{\Z}_2)$ comme \'etant pour chaque entier $k$ l'ensemble de niveaux $\lambda_k$ \`a partir duquel tout $k-$cycle ${\Z}_2$ dans $Z^0_{m-1,rel}(M^m, \p M^m; {\Z}_2)$ est homologue \`a un cycle dont la masse maximale est $\lambda_k$. La caract\'erisation du spectre non lin\'eaire, ou {\bf largeurs de Gromov}, au moyen du cup-produit est donn\'ee ainsi\footnote{Dans la litt\'erature anglophone sur le sujet on retrouve souvent \`a la place de  $\lambda_k$ la notation  $w_k$ pour ``width'' (``largeur'' en fran\c cais). Nous choisissons
$\lambda_k$ pour \^etre plus proche de l'ild\'ee originale de Gromov de g\'en\'eraliser la notion de spectre \`a des situations non lin\'eaires. La lettre grecque $\lambda$  raisonne aussi mieux avec le l.}
\be
\label{spec}
\lambda_k(M^m):=\inf_{\lf\{\Psi\, :\, X\rightarrow {\mathcal M}\ ;\ \Psi^\ast\al^k\ne 0 \ \mbox{dans } H^k(X,{\Z}_2) \rg\}}\quad\sup_{x\in X}\quad {\mathbf M}(\Psi(x))
\ee
o\`u $X$ varie librement sur les complexes simpliciaux de dimension finie  et $\al^k=\al\smile\cdots\smile\al$ ($k$ fois). Les applications $\Psi$  continues d'un complexe simplicial de dimension finie  dans ${\mathcal M}$ satisfaisant $\Psi^\ast\al^k\ne 0$ sont appel\'ees {\bf $k-$balayages}.  

Une premi\`ere remarque tautologique consiste \`a observer que $\Psi^\ast\al^k\ne 0$ donne pour tout entier $0<l<k$ $\Psi^\ast\al^l\ne 0$ et donc, ce que l'on appelle aussi le {\bf spectre volumique} est une suite croissante
\[
0<\la_1(M^m)\le \la_2(M^m)\cdots\le\la_k(M^m)\le\cdots
\]
Un exemple de $k-$balayage est
donn\'e par la restriction \`a ${\R}P^k$ de l'application $\Phi$ plus haut. On a en effet $\beta:=\Phi^\ast\al\ne 0$ et par ailleurs, un r\'esultat classique illustrant le cup produit sur les espaces projectifs r\'eels donne
 $\beta^k\ne 0$. Donc $\Phi^\ast \al^k\ne 0$.
 
 Dans le cas maintenant o\`u \underbar{$M^m$ est sans bord}, le th\'eor\`eme d'Almgren et Pitts s'\'etend\footnote{Voir le travail important de Marques et de Neves  \cite{MN} \`a ce sujet qui vient pr\'eciser et compl\'eter la th\'eorie.} aux $k-$balayages et on d\'emontre qu'en dimension plus petite que $8$ pour tout $k\in {\N}^\ast$ il existe une famille de surfaces plong\'ees disjointes 
 $\Sigma_1^k\cdots\Sigma_{n_k}^k$ et des multiplicit\'es enti\`eres non nulles $q^k_1\cdots q^k_{n_k}$ telles que
 \be
 \label{4a}
 \lambda_k(M^m,g)=\sum_{i=1}^{n_k}q_i^k\ |\Sigma_i^k|\quad.
 \ee
 La question principale de cet expos\'e devient alors celle de savoir si, lorsque $k$ tend vers l'infini,  la famille $\Sigma_i^k$ appartient \`a un ensemble fini ou infini de surfaces plong\'ees g\'eom\'etriquement distinctes.
 
En ayant \`a l'esprit les preuves de r\'esultats sur l'existence d'une infinit\'e de g\'eod\'esiques sur une vari\'et\'e donn\'ee, il est naturel d'\'etudier le comportement asymptotique du spectre $\la_k$.
Le r\'esultat suivant du principalement \`a Gromov et revisit\'e par Larry Guth dans \cite{Gut} est central pour la d\'emonstration du th\'eor\`eme ultime de notre expos\'e.

\begin{theo}\cite{Gro2}, \cite{Gro3}, \cite{Gut} 
\label{asympspec}
Soit $(M^m,g)$ une vari\'et\'e riemannienne sans bord donn\'ee alors lorsque $k$ tend vers l'infini on a le comportement asymptotique suivant\footnote{La diff\'erence de comportement asymptotique
du spectre si l'on compare (\ref{2}) avec (\ref{5}) dans le cas des g\'eod\'esiques, qui sont des \'el\'ements de $Z_{m-1}(M^m,{\Z}_2)$ pour $m=2$ pourrait surprendre \`a premi\`ere vue. Nous allons voir que le choix du groupe ${\Z}_2$ est responsable du fait que les changements de topologies dans $Z_1(M^2,{\Z}_2)$ arrivent ``plus vite'' au rythme de $\simeq\sqrt{k}$ pour la longueur alors que c'\'etait $\simeq k$ en multiplicit\'e enti\`ere.}
\be
\label{5}
\frac{\log\la_k(M^m)}{\log k}=\frac{1}{m}+O\lf(\frac{1}{\log k}\rg)\quad.
\ee 
\end{theo}
Les raisons de ce comportements asymptotiques m\'eritent d'\^etre un peu \'eclair\'ees \`a ce stade de notre pr\'esentation. 

La minoration  asymptotique $\la_k\ge C k^{1/m}$ vient d'un argument inspir\'e d'une preuve de Lusternik et Schnirelmann  (voir par exemple \cite{Kl} theorem 2.1.10). Soient $U_1\cdots U_k$ des ouverts disjoints 
de $M^m$ alors on d\'emontre l'in\'egalit\'e suivante
\be
\label{6}
\la_k(M^m)\ge\sum_{i=1}^k\la_1(U_i)\quad,
\ee
o\`u $\la_1(U_i)$ est la premi\`ere largeur de Gromov de l'ouvert $U_i$ dans $Z_{m-1}(U_i,\p U_i;{\Z}_2)$. En effet, soit $X$ un complexe simplicial de dimension fini et $\Psi$ de $X$ dans $M^m$ r\'ealisant un 
$k-$balayage : $\Psi^\ast \al^k\ne 0$ dans $H^k(X;{\Z}_2)$. On note
\[
X_i:=\lf\{x\in X\ ;\ {\mathbf M}(\Phi(x)\cap U_i)<\la_1(U_i)\rg\}
\]
Comme un balayage de $M^m$ r\'ealise un balayage de chaque sous ouvert de $M^m$, si $\iota_{X_i}\  :\ X_i\rightarrow X$ est l'inclusion canonique, on a donc $\iota_{X_i}^\ast \Psi^\ast\alpha=0$ dans $H^1(X_i,{\Z}_2)$. La suite exacte en cohomologie
\[
H^1(X,X_i;{\Z}_2)\xrightarrow{j_i^\ast} H^1(X;{\Z}_2)\xrightarrow{\iota_i^\ast} H^1(X_i,{\Z}_2)
\]
donne que $\al\in \mbox{Im}(j_i^\ast)$ et de plus il existe $\al_i\in H^1(X,X_i;{\Z}_2)$ tel que $j_i^\ast(\al_i)=\Psi^\ast\al$. Supposons par l'absurde que $X=\cup_{i=1}^kX_i$. On a donc trivialement
\[
H^k(X,\cup_{i=1}^k X_i;{\Z}_2)=0\quad.
\]
La propri\'et\'e fondamentale sur le cup produit relatif (page 209 \cite{Hat}) donne alors
\[
 H^1(X,X_1;{\Z}_2)\smile \cdots\smile H^1(X,X_k;{\Z}_2)\ \longrightarrow\ H^k(X,\cup_{i=1}^k X_i;{\Z}_2)=0
\]
donc en particulier $\al_1\smile\cdots\smile\al_k=0$. 

Pour tout $A, B\subset X$ le diagramme commutatif\footnote{$\Delta$ d\'esigne dans ce diagramme l'application diagonale $\Delta(x):=(x,x)$.}
\be
\label{comm}
\begin{array}{c}
\begin{array}{ccc}
X &\ds \xrightarrow{j_{A\cup B}} &(X,A\cup B)  \\[5mm]
\ds\downarrow{\Delta} &  &\ds \downarrow {\Delta}\\[5mm]
X\times X &\ds \xrightarrow{j_A\times j_B} &\quad(X,A)\times (X,B)= (X\times X, A\times X\cup X\times B) 
\end{array}
\end{array}
\ee
donne l'identit\'e $j^\ast_{A\cup B}\Delta^\ast=\Delta^\ast(j^\ast_A\times j_B^\ast)$ de laquelle on d\'eduit\footnote{ Au moyen de l'identit\'e fondamentale $a\smile b=\Delta^\ast (a\times b)$.}
\be
\label{cup}
j^\ast_{A\cup B}(a\smile b)=j_A^\ast a\smile j_B^\ast b\quad .
\ee
En it\'erant cette identit\'e pour $A:=\cup_{i=1}^lX_i$ et $B:=X_{l+1}$ on obtient finallement
\[
\Psi^\ast\al^k=j_{X_1}^\ast(\al_1)\smile\cdots j_{X_k}^\ast(\al_k)= j^\ast_{\cup X_i}(\al_1\smile\cdots\smile\al_k)=0\quad,
\]
ce qui est contradictoire. Ceci d\'emontre l'in\'egalit\'e (\ref{6}).

Une fois l'identit\'e (\ref{6}) \'etablie, on se convainc assez facilement que l'on peut  asymptotiquement d\'ecomposer la vari\'et\'e en $k$ cubes disjoints de diam\`etre $\simeq k^{-1/m}$ et tels que 
la largeur de Gromov de chacun de ces cubes soit environ $k^{(-m+1)/m}$. L'in\'egalit\'e (\ref{6}) donne alors $\la_k(M^m)\ge C\ k^{1/m}$, ce qui d\'emontre la minoration dans (\ref{5}).

Nous donnons une illustration de l'argument\footnote{Cette approche dite de ``bend and cancel'' est due \`a Larry Guth.}  dans le cadre le plus simple possible : $M^m:=T^2={\R}^2/{\Z}^2$. L'application
\[
\Psi_1 \ :\ X_1={\R}/{\Z}\ \longrightarrow\ Z_1(T^2;{\Z}_2)\quad \Psi(t):=\{x_1=t\}
\]
r\'ealise un balayage de Birkhoff. Plus g\'en\'eralement, on consid\`ere l'application suivante. Soit $\pi$ la projection qui \`a un point de la droite r\'eelle associe le point le plus proche dans $[0,1]$. On appelle alors
$\Psi_k$ l'application qui \`a un polyn\^ome non nul $P$ de degr\'e au plus $k$ associe l\'el\'ement de $Z_1({\R}^2/{\Z}^2;{\Z}_2)$ \'egal \`a la somme des cercles \`a multiplicit\'e ${\Z}_2$ donn\'es par $x_1= \pi(a_i)$
o\`u $a_i$ sont les racines r\'eelles de $P$.  On d\'emontre (voir \cite{Gut1} section 5) que $\Psi_k$ est continue\footnote{La continuit\'e est principalement due \`a la multiplicit\'e ${\Z}_2$ : deux racines r\'eelles qui convergent l'une vers l'autre s'annihilent \`a la limite.}
de ${\R}P^k$ \footnote{L'espace $X$ des polyn\^omes non nuls de deg\'e au plus $k$ quotient\'e par ${\R}^\ast$ est diff\'eomorphe \`a ${\R}P^k$.} dans $Z_1(T^2;{\Z}_2)$. L'application $\Psi_k$ envoie un g\'en\'erateur
du $H_1({\R}P^k,{\Z}_2)$ donn\'e par la suite de polyn\^ome $P_t(Y)= (t+t^{-1}) Y-2$ sur un g\'en\'erateur du $H_1(Z_1(T^2;{\Z}_2))$ et donc on a $\Psi_k^\ast\al\ne 0$. Ceci implique que $(\Psi_k^\ast\al)^k\ne 0$.
Or $(\Psi_k^\ast\al)^k=\Psi_k^\ast(\al^k)$, donc $\Psi_k$ r\'ealise un $k-$balayage de $Z_1(T^2;{\Z}_2)$.

Chaque \'el\'ement $\Psi_k(x)$ est constitu\'e au plus de $k$ cercles de longueur $1$. On a donc
\[
\max_{x\in {\R}P^k}{\mathbf M}(\Psi_k(x))=k\quad.
\]
Or ce n'est pas la borne optimale donn\'ee par le th\'eor\`eme qui est asymptotiquement $\simeq \sqrt{k}$ dans ce cas. C'est que nous n'avons pas fait l'usage complet du fait que nous avons une multiplicit\'e ${\Z}_2$ qui est la vrai ``responsable'' de cette croissance sous-lin\'eaire qui sera si utile par la suite. 

L'id\'ee de Larry Guth est de choisir un maillage de $T^2$ fait de $k$ cubes de taille identique $\simeq \sqrt{k}$ et de ``projeter'' $\Psi_k$ le maximum sur ce maillage. Il d\'emontre que chaque cercle de $\Psi_k(x)$ se projette int\'egralement sur ce maillage \`a l'exception d'une portion de cercle contenue dans exactement 2 cubes  et de taille $\simeq 1/\sqrt{k}$ et cela continument en $x$ (fig. 2). 

\psfrag{0}{$0$}
\psfrag{1}{$1$}
\psfrag{z}{z\'eros de $x$}
\psfrag{RZ}{${\R}^2/{\Z}^2$}
\psfrag{Mult}{Multiplicit\'e ${\Z}_2$}
\psfrag{psi}{$\psi_3(x)$}
\psfrag{x}{$x \rightarrow x + \delta x$}
\psfrag{ay}{$a_3 \,Y^3 + a_2 \,Y^2 + a_1 \,Y + a_0 / {\R}^*_+ = x \in \underset{\stackrel{^\parallel\;\,}{\!\!X_3}}{{\R} {\mathbb P}^3}$}
\begin{center}
\includegraphics[width=10cm]{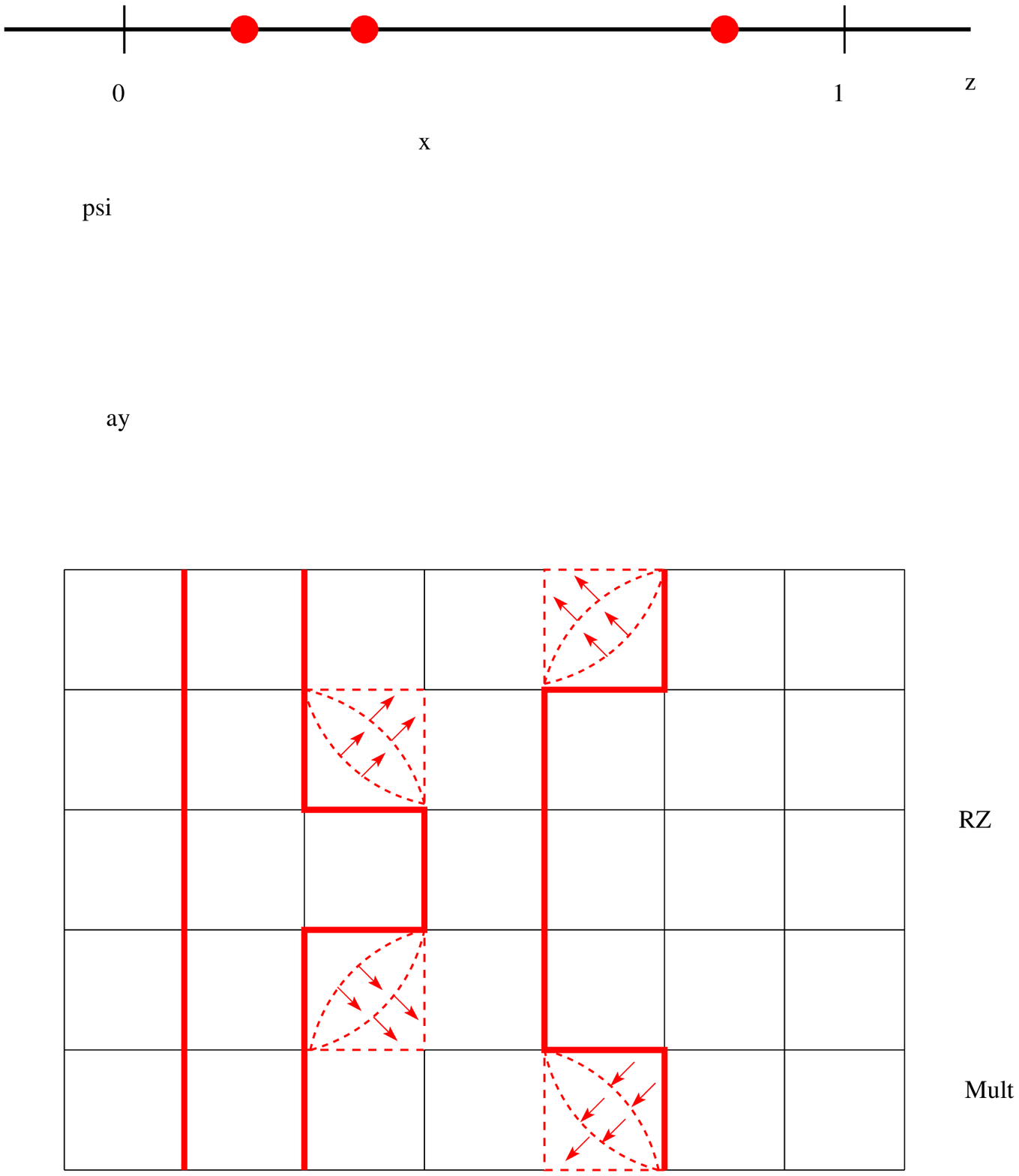}
\end{center}

\bigskip
\begin{center}
Fig.~2: ~ La proc\'edure ``plier et annuler'' de Guth
\end{center}

On obtient ainsi une nouvelle application $\tilde{\Psi}_k$ v\'erifiant
\[
\max_{x\in {\R}P^k}{\mathbf M}(\ti{\Psi}_k(x))\le C\  k\times \frac{1}{\sqrt{k}}+ {\mathbf M}(\mbox{maillage})\le \ C'\ \sqrt{k}\quad,
\]
ce qui donne l'estimation optimale.

\section{La loi de Weyl de Liokumovich,  Marques et Neves.}
R\'ecemment Yevgeny Liokumovich, Fernando Cod\'a Marques et Andr\'e Neves sont parvenus \`a pr\'eciser l'estimation (\ref{5}) et \`a  \'etablir une {\bf loi de Weyl}\footnote{La loi de Weyl pour le {\bf spectre de l'op\'erateur de Laplace Beltrami}  sur une vari\'et\'e riemanniene ferm\'ee $M^m$ donne l'existence d'une constante $c(m):=4\pi^2\ |B^m_1|^{-2/m}$ telle que
\[
\lim_{k\rightarrow +\infty}\la_k\ k^{-2/m}=c(m) \ \mbox{vol}(M^m)^{-2/m}\quad.
\] } pour le {\bf spectre du volume} tout \`a fait remarquable\footnote{Le th\'eor\`eme-\ref{loiweyl} a r\'epondu positivement \`a une conjecture de Gromov.}
\begin{theo}
\label{loiweyl}
Soit $m$ un entier sup\'erieur \`a 1. Il existe une constante universelle\footnote{Une expression explicite de cette constante m\^eme dans le cas de la plus basse dimension $m=2$  au moyen de quantit\'es ``connues'' est un probl\`eme ouvert.}  $a(m)>0$ telle que pour toute vari\'et\'e compacte riemannienne $(M^m,g)$ (avec ou sans bord) on ait
\be
\label{7}
\lim_{k\rightarrow +\infty}\la_k\ k^{-1/m}=a(m) \ \mbox{vol}(M^m)^{(m-1)/m}\quad.
\ee 
o\`u $\lambda_k$ est le spectre du volume donn\'e par (\ref{spec}).
\end{theo}
Au del\`a de l'aspect incontestablement esth\'etique de ce r\'esultat, la pr\'ecision de celui-ci induit un m\'ecanisme qui va g\'en\'erer de nombreuses surfaces minimales r\'epondant ainsi partiellement \`a une conjecture ancienne de Shing-Tung Yau.
\begin{conj}
\label{yau}\cite{Yau} Toute vari\'et\'e ferm\'ee 3-dimensionnelle poss\`ede un nombre infini de surfaces minimales immerg\'ees g\'eom\'etriquement distinctes. 
\end{conj}
Pr\'ecisemment dans un article en collaboration avec Kei Irie, Fernando Cod\'a Marques et Andr\'e Neves d\'emontrent au moyen de la loi de Weyl (\ref{7}) le th\'eor\`eme suivant.
\begin{theo}\cite{IMN}
\label{infweyl} Soit $M^m$ une vari\'et\'e ferm\'ee de dimension $3\le m\le 7$. Alors pour un ensemble g\'en\'erique de m\'etriques l'union des hypersurfaces minimales ferm\'ees et plong\'ees 
est dense dans $M^m$.
\end{theo}
La strat\'egie de la preuve peut-\^etre comprise ainsi. Soit $U$ un ouvert de $M^m$. On note ${\mathfrak M}_U$ l'espace des m\'etriques pour lesquelles il existe une surface minimale plong\'ee
non d\'eg\'en\'er\'ee\footnote{La d\'eriv\'ee seconde de l'aire n'a pas de noyau.} intersectant $U$. On d\'emontre sans trop de difficult\'es, en utilisant le th\'eor\`eme des fonctions implicites, combin\'e aux r\'esultats d\'esormais classiques
de White \cite{Whi} sur l'espace de surfaces minimales pour les m\'etriques g\'en\'eriques, que ${\mathfrak M}_U$ est ouvert dans l'espace ${\mathfrak M}$ des m\'etriques $C^\infty$.

Le point central de la d\'emonstration consiste \`a \'etablir le fait que ${\mathfrak M}_U$ est dense dans ${\mathfrak M}$. Supposons que ce ne soit pas le cas et prenons une m\'etrique g\'en\'erique $g$ qui poss\`ede
un voisinage ouvert ${\mathfrak V}$ dans ${\mathfrak M}$ pour lequel il n'existe pas de surface minimale plong\'ee non d\'eg\'en\'er\'ee intersectant $U$ pour toute m\'etrique dans ${\mathfrak V}$. Pour une m\'etrique g\'en\'erique un r\'esultat de Ben Sharp \cite{Sha} affirme que l'ensemble des surfaces minimales non d\'eg\'en\'er\'ees et affect\'ees de multiplicit\'e enti\`eres quelconques r\'ealise un ensemble d\'enombrable. 
On consid\`ere alors une perturbation de $g$ de la forme $g_t:=(1+t\,\chi)\, g$ pour $t>0$ et $\chi$ est une fonction quelconque non n\'egative et non nulle support\'ee dans $U$. Pour $t$ petit $g_t\in {\mathfrak V}$. 
Supposons que $(M^m,g_t)$ ne poss\`ede aucune surface minimale plong\'ee et ferm\'ee intersectant $U$. Cela signifie que toutes les valeurs de $\la_k(M^m,g_t)$ sont r\'ealis\'ees par des surfaces minimales de $g$ et donc
 l'application $t\rightarrow \la_k(M^m,g_t)$ pour $t$ petit ne peut prendre qu'un nombre au plus d\'enombrables de valeurs. Or on d\'emontre sans trop d'efforts que celle-ci est continue. Cela implique alors qu'elle est constante proche de zero. La loi de Weyl (\ref{7}) donne cependant pour $k$ suffisamment grand
 \[
\lambda_k(M^m,g_t)>\lambda_k(M^m,g)\quad,
\]
ce qui est contradictoire. Donc $g_t$ poss\`ede une surface minimale $\Sigma$ plong\'ee et ferm\'ee passant par $U$. Un r\'esultat qui remonte \`a W.Klingenberg dans le cas des g\'eod\'esiques (voir \cite{Kl} capitre 3) affirme que l'on peut approcher la m\'etrique $g_t$ par une suite de m\'etriques (en la modifiant \`a l'ordre 2 et plus au voisinage de $\Sigma$) de telle sorte que pour chacune d'entre elles $\Sigma$ soit toujours une surface minimale mais non d\'eg\'en\'er\'ee cette fois. On a donc contredit le fait que ${\mathfrak V}$ existe et
on en d\'eduit que ${\mathfrak M}_U$ est dense dans ${\mathfrak M}$ pour tout ouvert $U$ de $M^m$.

Soit $(U_i)_{i\in {\N}}$ une base d\'enombrable de la topologie de $M^m$. Le sous espace de ${\mathfrak M}$ donn\'e par $\cap_{i\in {\N}}{\mathfrak M}_{U_i}$ est solution du probl\`eme grace au th\'eor\`eme de Baire. Ceci  conclut la preuve du th\'eor\`eme~\ref{infweyl}.

\section{La propri\'et\'e de Frankel}

Un th\'eor\`eme d\^u \`a Jacques Hadamard \cite{Had} affirme que sur une surface ferm\'ee \`a courbure strictement positive les deux \'el\'ements de toutes paires arbitraires de g\'eod\'esiques ferm\'ees s'intersectent\footnote{En fait dans \cite{Had} il est d\'emontr\'e que, sous les m\^emes hypoth\`eses, toute g\'eod\'esique ferm\'ee est intersect\'ee un nombre infini de fois par toute autre g\'eod\'esique compl\`ete mais non ferm\'ee}. 

Ce r\'esultat a \'et\'e g\'en\'eralis\'e au cas des hypersurfaces par Th\'eodore Frankel.
\begin{theo}
\label{frank}\cite{Fra}
Soit $(M^m,g)$ une vari\'et\'e riemannienne compl\`ete et connexe dont la courbure de Ricci est strictement positive. Alors deux hypersurfaces minimales compl\`etes arbitraires dont l'une est compacte doivent s'intersecter.
\end{theo}
La d\'emonstration se comprend ainsi. On suppose pour la simplicit\'e de la pr\'esentation que les deux surfaces minimales $\Sigma_0$ et $\Sigma_1$ sont plong\'ees. Soit $\gamma\ :\ [0,1]\ \rightarrow\ M^m$ une g\'eod\'esique 
joignant les deux surfaces et minimisant la distance entre elles. Un argument simple montre que $\gamma$ doit inters\'ecter $\Sigma_0$ et $\Sigma_1$ chacune orthogonalement en deux points $p_0$ et $p_1$.
La g\'eod\'esique est stable et un calcul direct de la d\'eriv\'ee seconde de la longueur en $\gamma$ donne pour une perturbation $e(t)\in \Gamma(\gamma^\ast TM^m)$ obtenue en transportant parall\`element un vecteur $e(0)\in T_{p_1}\Sigma_1$ de $p_0$ le long de $\gamma$ jusqu'\`a $p_1$
\be
\label{8}
0\le D^2L_\gamma( e,e)={\mathbb I}^{\,1}_{\,p_1}(e(1),e(1))-{\mathbb I}^{\,0}_{\,p_0}(e(0),e(0))-\int_0^1 K_{\gamma(t)}(e(t),\dot{\gamma}(t)) \ |\dot{\gamma}(t)|^{-1}\ dt
\ee
o\`u ${\mathbb I}^{j}_{p_j}$ est la seconde forme fondamentale de la surface $\Sigma_j$ au point $p_j$ (pour $j=0,1$)dans $(M^m,g)$ et $K_x(X,Y)$ est la courbure sectionnelle de $M^m$ au point $x$ suivant le couple de vecteur $(X,Y)$ tangents en ce point. 
Prenons maintenant non seulement un seul vecteur $e(0)$ mais une base orthonormale de $T_{p_0}\Sigma_0$ : $(e_i(0))_{i=1\cdots m-1}$. Cette base, transport\'ee paral\`ellement et \`a laquelle on adjoint le vecteur unit\'e tangent de $\gamma$, r\'ealise un rep\`ere mobile orthonorm\'e de $\gamma^\ast TM^m$ (fig. 3). 

\begin{center}
\psfrag{S0}{$\Sigma_0$}
\psfrag{S1}{$\Sigma_1$}
\psfrag{e1(0)}{$e_1(0)$}
\psfrag{e2(0)}{$e_2(0)$}
\psfrag{e1(1)}{$e_1(1)$}
\psfrag{e2(1)}{$e_2(1)$}
\psfrag{g}{$\gamma$}
\includegraphics[width=8cm]{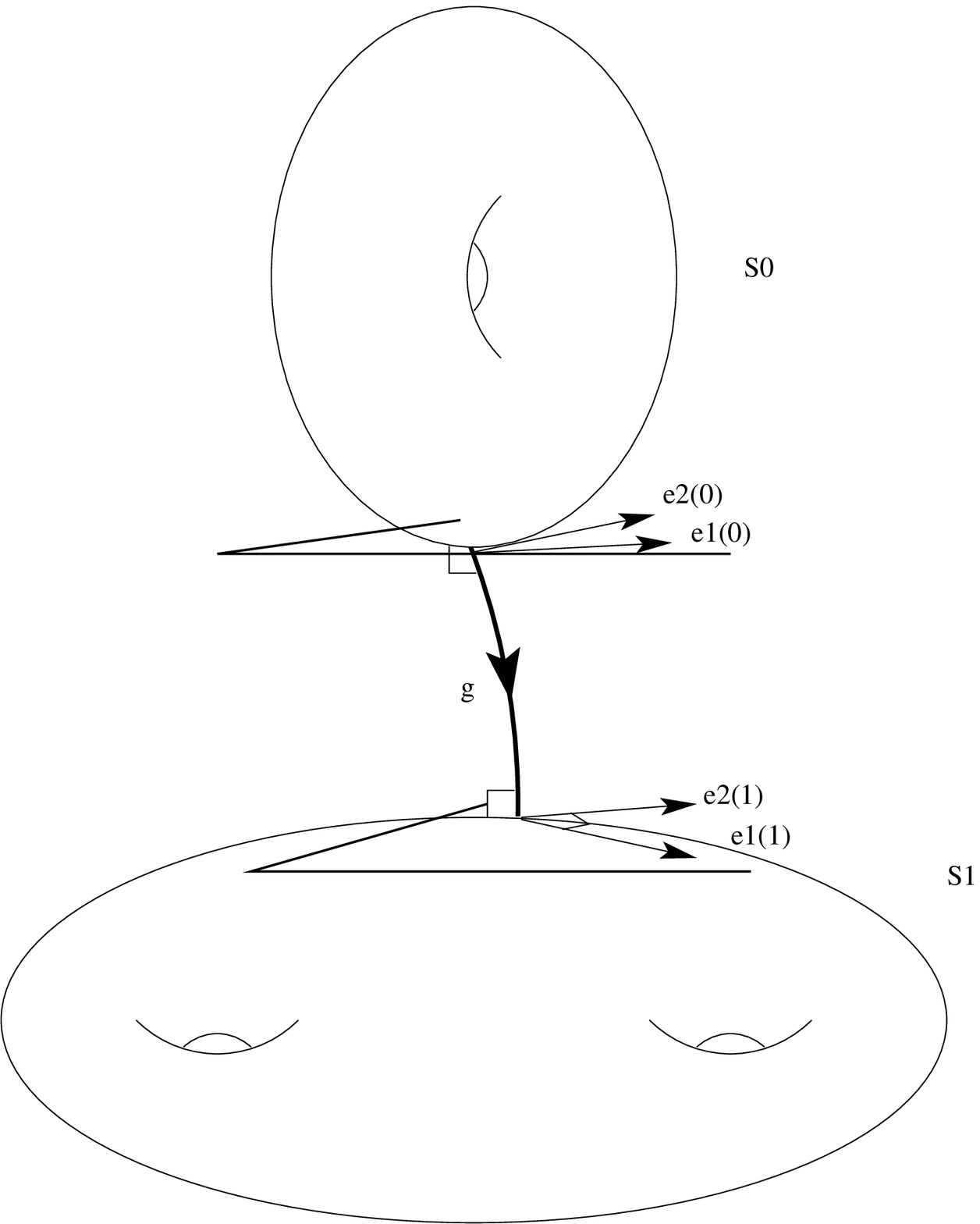}
\end{center}

\bigskip
\begin{center}
Fig.~3: ~ La propri\'et\'e de Frankel
\end{center}
 En particulier, 
\[
\sum_{i=1}^{m-1}K_{\gamma(t)}(e_i(t),\dot{\gamma}) =Ric_{\gamma(t)}(\dot{\gamma})
\]

Comme $\Sigma_j$ est une surface minimale on a par ailleurs l'annulation de la courbure moyenne en $p_j$ ce qui donne 
\[
\sum_{i=1}^{m-1}{\mathbb I}^{\,0}_{\,p_0}(e_i(0),\dot{\gamma}) =0\quad\mbox{ et }\quad \sum_{i=1}^{m-1}{\mathbb I}^{\,1}_{\,p_1}(e_i(1),\dot{\gamma}) =0\quad.
\]
En additionnant (\ref{8}) pour chaque vecteur $e_i$ on obtient
\[
0\le \sum_{i=1}^{m-1} D^2L_\gamma( e_i,e_i)=-\int_0^1\ Ric_{\gamma(t)}(\dot{\gamma})\ |\dot{\gamma}(t)|^{-1}\ dt
\]
L'hypoth\`ese de courbure de Ricci strictement positive contredit cette in\'egalit\'e et on conclut que les deux hypersurfaces s'intersectent n\'ecessairement.

Inspir\'e par ce r\'esultat, on d\'efinit la  {\bf propri\'et\'e de Frankel} pour une vari\'et\'e riemannienne $M^m$ comme \'etant la propri\'et\'e suivant laquelle toute paire de surfaces minimales ferm\'ees doivent avoir une intersection non vide.

 Marques et Neves d\'emontrent le r\'esultat suivant.
\begin{theo}
\label{frankinf}\cite{MN}
Soit $(M^m,g)$ une vari\'et\'e riemannienne compacte  de dimension $3\le m\le 7$ satisfaisant \`a la propri\'et\'e de Frankel. Alors $M^m$ poss\`ede une infinit\'e d'hypersurfaces minimales plong\'ees ferm\'ees et g\'eom\'etriquement distinctes.
\end{theo}
En combinant les deux th\'eor\`emes pr\'ec\'edents on obtient en particulier le corollaire suivant.
\begin{coro}
\label{riccipos}
Soit $(M^m,g)$ une vari\'et\'e riemannienne compacte  de dimension $3\le m\le 7$ et de courbure de Ricci strictement positive, alors $M^m$ poss\`ede une infinit\'e d'hypersurfaces minimales plong\'ees ferm\'ees et g\'eom\'etriquement distinctes.
\end{coro}
La d\'emonstration du th\'eor\`eme~\ref{frankinf} peut se comprendre ainsi. La propri\'et\'e de Frankel et le fait que les diff\'erentes surfaces $\Sigma^k_i$ dans (\ref{4a}) contribuant \`a la r\'ealisation de la valeur propre volumique $\la_k$ sont disjointes impose $n_k=1$. Donc pour tout $k\in {\N}^\ast$ il existe un entier $q^k\in {\N}^\ast$ et une surface plong\'ee $\Sigma^k$ tels que
\[
\la_k=q^k\, |\Sigma^k|\quad.
\]
Supposons dans un premier temps que les valeurs propres volumiques soient toutes diff\'erentes et donc $0<\la_1<\la_2<\cdots<\la_k\cdots$. Si il n'existait dans $M^m$ qu'un nombre fini $N$ d'hypersurfaces minimales plong\'ees dont l'aire minimale est $A>0$ alors $\la_k\ge k A/N$ ce qui contredit la croissance sous lin\'eaire\footnote{ C'est ici que le choix de la multiplicit\'e ${\Z}_2$ pour les cycles joue \`a plein son r\^ole. Nous avons expliqu\'e plus haut comment ce choix etait \`a l'origine de la croissance sous-lin\'eaire du spectre.} du spectre donn\'ee par (\ref{5}).

Il reste \`a consid\'erer le cas o\`u deux valeurs propres succ\'essives sont \'egales : $\la_k=\la_{k+1}$ pour un $k\in {\N}^\ast$. 

Ce cas est un cas bien connu de la th\'eorie\footnote{Il s'agissait alors dans cette th\'eorie  du spectre g\'en\'er\'e par les changements d'homologie des ensembles de niveau ${\mathcal M}^\la$ o\`u ${\mathcal M}=W^{1,2}(S^1,M^m)$ comme nous l'avons vu dans la premi\`ere partie de l'expos\'e. }de Liusternik-Shnirelmann qui permet d'affirmer l'existence d'une infinit\'e de g\'eod\'esiques ferm\'ees. L'argument reprend les propri\'et\'es fondamentales du cup-produit dans l'esprit  de la d\'emonstration de l'in\'egalit\'e (\ref{6}).

 Plus  pr\'ecisemment, supposons que l'espace des hypersurfaces minimales plong\'ees se r\'eduit \`a un nombre fini de surfaces $\Sigma_1\cdots \Sigma_N$. 
 
 On d\'emontre dans un premier temps sans trop de difficult\'es qu'il existe $\ep>0$ tel qu'aucune application $\Xi$ d'un complexe simplicial de dimension fini arbitraire $Z$ \`a valeur dans l'union des boules $B_\ep^{\mathcal F}(\Sigma_i)$ ne r\'ealise un balayage de $M^m$.  Cela signifie en d'autre termes que pour un tel $\Xi$ on a n\'ec\'essairement $\Xi^\ast\al=0$.
 
  \`A tout $(k+1)-$balayage $\Psi\ :\ X\longrightarrow {\mathcal M}:=Z^0_{m-1}(M^m,{\Z}_2)$ on associe l'ensemble $$Y:=\Psi^{-1}\lf({\mathcal M}\setminus \cup_{i=1}^NB_\ep^{\mathcal F}(\Sigma_i)\rg)$$
On d\'emontre dans un premier temps que $\Psi\ :\ Y\rightarrow {\mathcal M}$ r\'ealise un $k-$balayage. En effet, si tel n'\'etait pas le cas, on aurait $\iota_Y^\ast(\Phi^\ast\al^k))=0$ o\`u $\iota_Y$ est l'inclusion canonique de $Y$ dans $X$.
La suite exacte en cohomologie
\[
H^k(X,Y;{\Z}_2)\xrightarrow{j_Y^\ast} H^k(X;{\Z}_2)\xrightarrow{\iota_Y^\ast} H^k(Y;{\Z}_2)
\]
donne alors l'existence de $\beta\in H^k(X,Y;{\Z}_2)$ tel que $j_Y^\ast\beta=\Phi^\ast\al^k$. Le nombre $\epsilon>0$ a \'et\'e choisi de fa\c con \`a ce que $\iota_{X\setminus Y}^\ast\Phi^\ast\al=0$. La suite exacte
de cohomologie suivante
\[
H^1(X,X\setminus Y;{\Z}_2)\xrightarrow{j_{X\setminus Y}^\ast} H^1(X;{\Z}_2)\xrightarrow{\iota_{X\setminus Y}^\ast} H^1(X\setminus Y;{\Z}_2)
\]
donne l'existence de $\gamma\in H^1(X,X\setminus Y;{\Z}_2)$ tel que $j_{X\setminus Y}^\ast\gamma=\Phi^\ast\al$. La propri\'et\'e fondamentale sur le cup produit relatif donne alors
\[
\beta\smile \gamma\in H^{k+1}(X,X;{\Z}_2)=0
\]
Si on applique (\ref{cup}) \`a $A=Y$, $B=X\setminus Y$, $a= \beta$ et $b=\gamma$ on obtient
\[
0=j^\ast_{X}(\beta\smile \gamma)= j_Y^\ast\beta\smile j_{X\setminus Y}^\ast\gamma=\Phi^\ast \al^{k+1}
\]
ce qui contredit le fait que $\Phi$ est un $k+1-$balayage. On en d\'eduit que la restriction de $\Phi$ \`a $Y$ est un $k-$balayage et donc par hypoth\`ese
\[
\sup_{y\in Y}\Phi(y)\ge \la_k=\la_{k+1}\quad.
\]
Les arguments variationnels de type minmax de la th\'eorie de Almgren-Pitts-Marques-Neves  construits \`a partir des restrictions de $\Phi$ aux sous simplexes $Y$ (voir \cite{MN} page 605) permettent alors de mettre en \'evidence une surface minimale plong\'ee en dehors des boules $B_\ep^{\mathcal F}(\Sigma_i)$, ce qui est une contradiction.
\section{La preuve de Song.}

La derni\`ere partie de cet expos\'e est consacr\'ee \`a la d\'emonstration de la conjecture de Yau  (dans sa version ``forte'' avec des plongements et non pas des immersions) en dimension inf\'erieure ou \'egale \`a $7$ par Antoine Song. Pr\'ecisemment
nous allons d\'emontrer le r\'esultat suivant.

\begin{theo}
\label{th-song}\cite{Son}
Soit $(M^m,g)$ une vari\'et\'e riemannienne ferm\'ee quelconque de dimension $3\le m\le 7$, alors $(M^m,g)$ poss\`ede une infinit\'e d'hypersurfaces minimales plong\'ees ferm\'ees et g\'eom\'etriquement distinctes.\hfill $\Box$
\end{theo}

Le reste de cette section est consacr\'e \`a la description de la preuve donn\'ee par Song.

Si la propri\'et\'e de Frankel est v\'erifi\'ee alors le th\'eor\`eme~\ref{frankinf} donne le r\'esultat. Nous supposons donc que la propri\'et\'e de Frankel n'est pas v\'erifi\'ee et, par l'absurde, que $M^m$ ne contient qu'un nombre fini d'hypersurfaces minimales plong\'ees. 

L'id\'ee principale consiste dans un premier temps \`a ``extraire'' de
$M^m$ un ouvert $U$ non vide, appel\'e le {\bf coeur}, bord\'e par des hypersurfaces minimales localement minimisantes dans $U$,  et dans lequel les deux \'el\'ements de toute paire de surfaces minimales plong\'ees, et sans bord dans $\ov{U}$, doivent s'intersecter. Cet ouvert $U$ va \^etre obtenu par une proc\'edure r\'ecurrente suivante.

Pour simplifier la pr\'esentation nous supposons que les hypersurfaces minimales plong\'ees de $M^m$ (qui sont en nombre fini par hypoth\`ese)  ont toutes deux faces\footnote{On peut toujours se ramener \`a cette situation en travaillant avec un rev\`etement double de $M^m$.}  dans $M^m$ et sont donc orientables. 

On d\'emontre au moyen du th\'eor\`eme des fonctions implicites\footnote{Un r\'esultat s'approchant de cette affirmation est explicitement d\'emontr\'e dans \cite{Koi}.} que chaque hypersurface minimale plong\'ee et ferm\'ee de  $M^m$ appartient \`a l'une des trois cat\'egories suivantes (fig. 4) :
\begin{itemize}
\item[-] Surfaces \`a voisinage {\bf contractant}. Ce sont les surfaces $S$ qui poss\`edent un voisinage \'equip\'e d'un feuillage r\'egulier de codimension 1 contenant la surface et auquel est associ\'e un champs de vecteur courbure moyenne $\vec{H}$ pointant\footnote{Il est important de rappeler \`a ce stade que 
\[
\delta\mbox{Aire}_S \,\vec{w}:=-\int_S \vec{H}\cdot \vec{w}\ dvol_{S}
\] C'est \`a dire que le champs de vecteur moyenne sur une surface est l'oppos\'e du ``gradient'' de l'aire. Traditionellement, en th\'eorie de mesure g\'eom\'etrique, la courbure moyenne est la somme des courbures principales  sans que cette somme soit divis\'ee par le nombre $m-1$ de directions. Cela rompt malheureusement avec la convention des g\'eom\`etres introduite dans l'article de Jean-Baptiste Meusnier de 1776 ! } vers $S$.
\medskip
\item[-] Surfaces \`a voisinage {\bf dilatant}. Ce sont les surfaces $S$ qui poss\`edent un voisinage \'equip\'e d'un feuillage r\'egulier de cdimension 1 contenant la surface et auquel est associ\'e un champs de vecteur courbure moyenne  pointant vers l'ext\'erieur de $S$.
\medskip
\item[-] Surfaces \`a voisinage {\bf mixte}. Ce sont les surfaces $S$ orientables qui poss\`edent un voisinage \'equip\'e d'un feuillage r\'egulier de codimension 1 contenant la surface et auquel est associ\'e un champs de vecteur courbure moyenne pointant vers l'ext\'erieur de $S$ d'un c\^ot\'e de $S$ et pointant vers $S$ de l'autre c\^ot\'e.
\end{itemize}
\begin{center}
\psfrag{S1}{$\Sigma_1$}
\psfrag{S2}{$\Sigma_2$}
\psfrag{S3}{$\Sigma_3$}
\psfrag{S4}{$\Sigma_4$}
\psfrag{d}{dilatant}
\psfrag{c}{contractant}
\psfrag{m}{mixte}
\includegraphics[width=12cm]{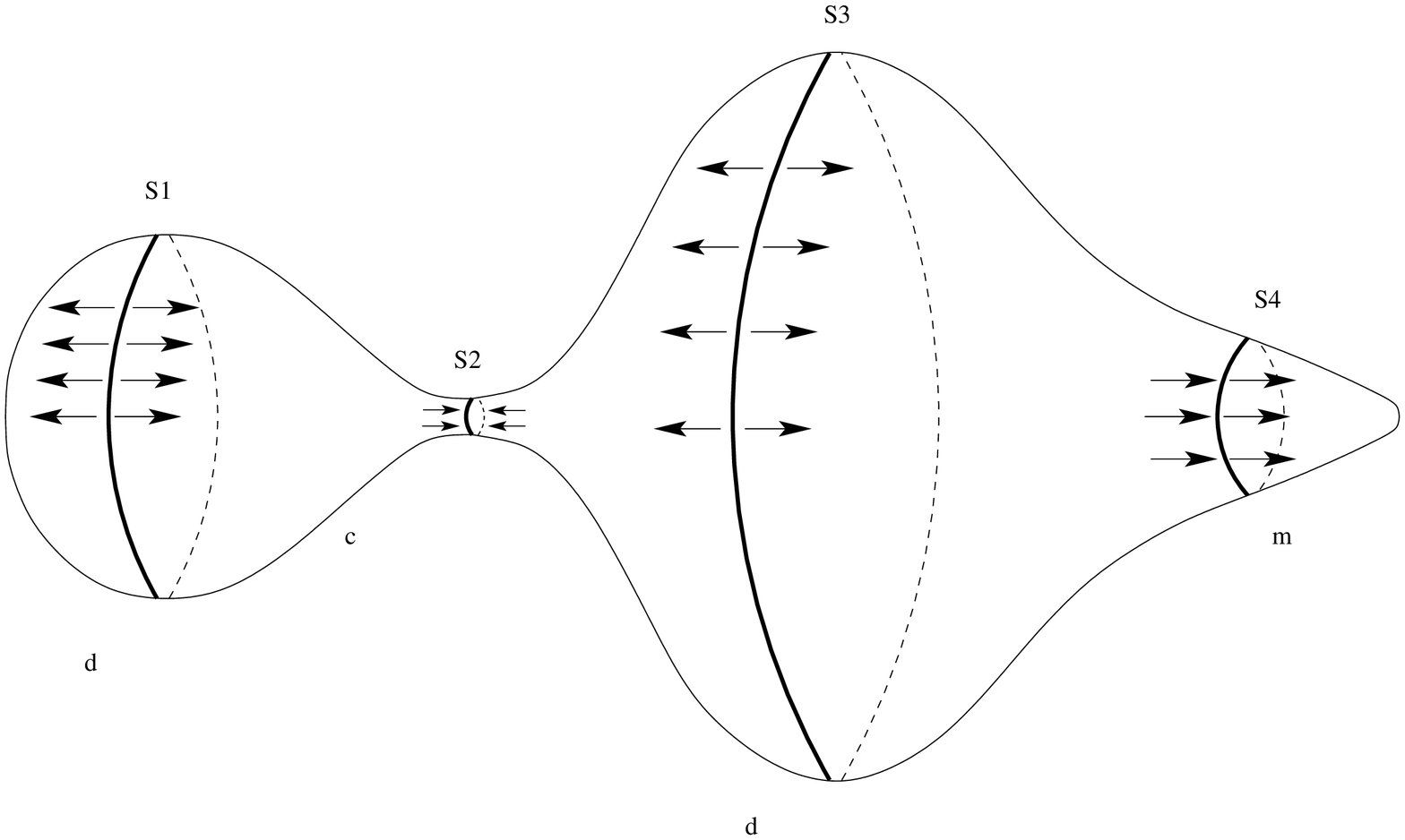}
\end{center}

\bigskip
\begin{center}
Fig.~4: ~ Les 3 familles de surfaces minimales
\end{center}

\begin{rema}Il est important pour la suite de remarquer \`a ce stade que toute surface minimale plong\'ee ferm\'ee \`a voisinage dilatant ou mixte ne peut \^etre localement minimisante pour l'aire.
\end{rema}
La construction par r\'ecurrence du coeur\footnote{Le coeur \`a priori n'est pas unique et va d\'ependre de choix faits au fur et \`a mesure de la proc\'edure.} d\'emarre ainsi. Soient $\Sigma_1$ et $\Sigma_2$ deux surfaces minimales plong\'ees et ferm\'ees de $M^m$ ne s'intersectant pas.  Il se pr\'esente alors quatre possibilit\'es

{\bf 1\`ere possibilit\'e.} Si l'une des deux surfaces, par exemple $\Sigma_1$, est \`a voisinage {\bf contractant} alors on d\'ecoupe $M^m$ suivant $\Sigma_1$. Soit $U_1$ une composante\footnote{Il peut y en avoir une ou deux.} connexe arbitraire de $M^m\setminus \Sigma_1$.

\begin{center}
\psfrag{S1}{$\Sigma_1$}
\psfrag{S2}{$\Sigma_2$}
\psfrag{c}{contractant}
\psfrag{U1}{$U_1$}
\includegraphics[width=8cm]{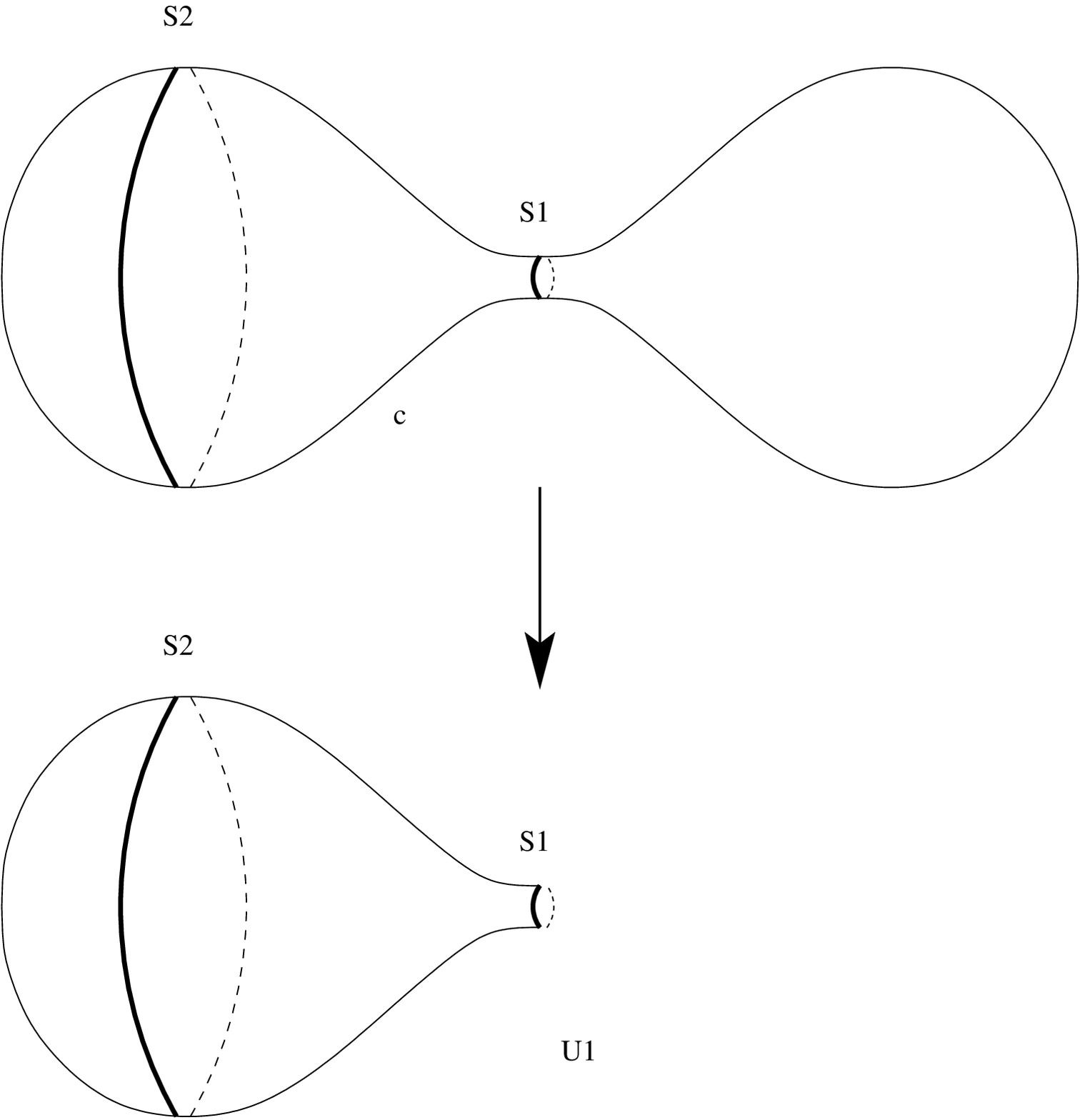}
\end{center}

\bigskip
\begin{center}
Fig.~5: ~ La construction du coeur : 1\`ere possibilit\'e.

\end{center}

{\bf 2\`eme possibilit\'e.} Si une des deux surfaces, par exemple $\Sigma_1$, est \`a voisinage {\bf mixte et si $M^m\setminus \Sigma_1$ a deux composantes connexes} alors on choisit $U_1$ comme \'etant la composante contractante c'est \`a dire qui contient le voisinage de $\Sigma_1$ avec le feuilletage dont le champs de courbure moyenne associ\'e pointe vers $\Sigma_1$.

\begin{center}
\psfrag{S1m}{$\Sigma_1$ mixte}
\psfrag{S2}{$\Sigma_2$}
\psfrag{c}{contractant}
\psfrag{U1}{$U_1$}
\includegraphics[width=5.5cm]{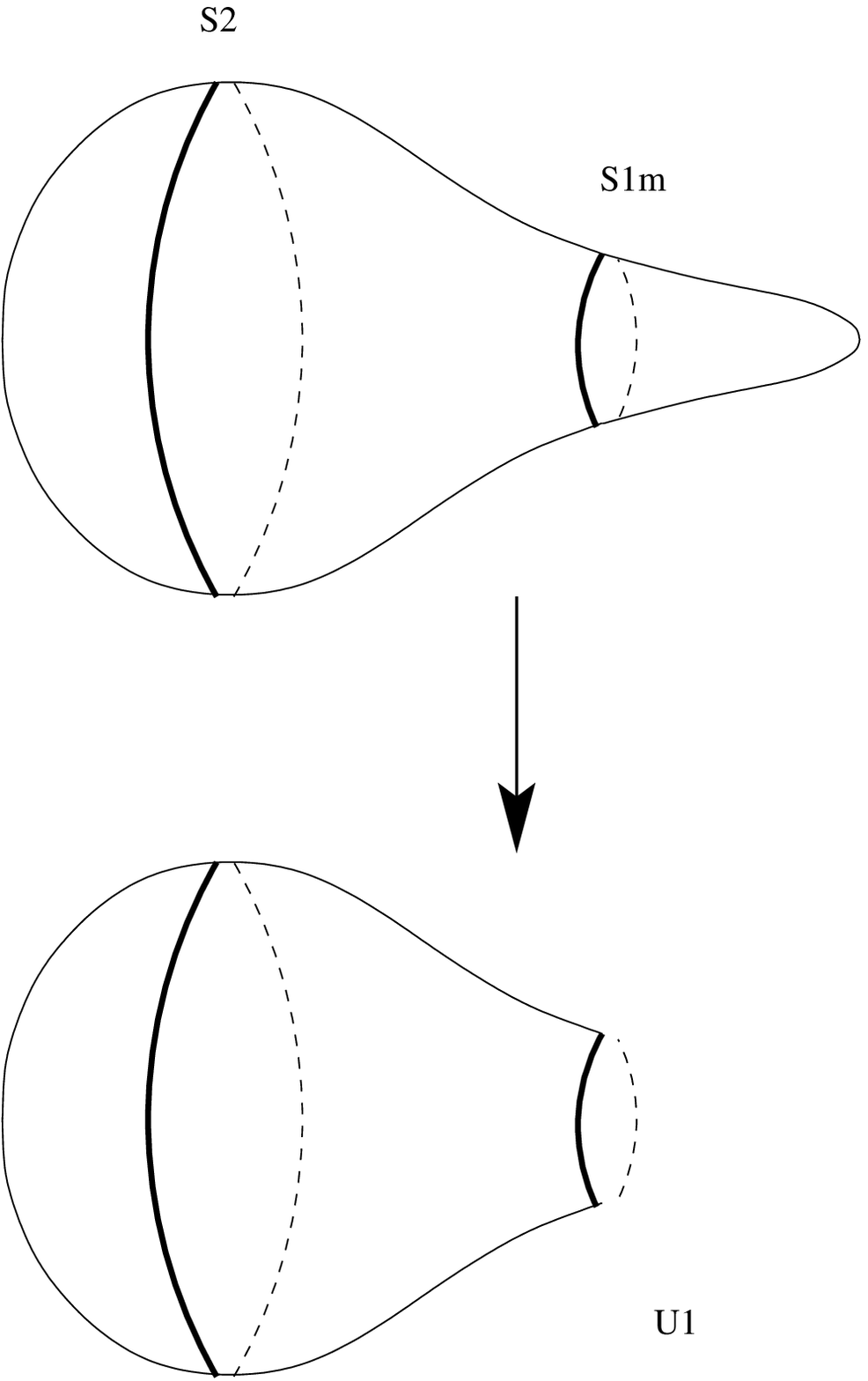}
\end{center}

\bigskip
\begin{center}
Fig.~6: ~ La construction du coeur : 2\`eme possibilit\'e.
\end{center}

{\bf 3\`eme possibilit\'e.} Si une des deux surfaces, par exemple $\Sigma_1$ n'est {\bf pas ``s\'eparante''} c'est \`a dire que  $M^m\setminus \Sigma_1$ a une seule composante connexe alors la classe dans $H_{m-1}(M^m;{\Z})$ d\'efinit par $\Sigma_1$ est non nulle. La minimisation de l'aire dans cette classe et le th\'eor\`eme~\ref{reg} donne une surface minimale $S$ plong\'ee ferm\'ee et n\'ecessairement \`a voisinage contractant. On choisit alors $U_1:=M^m\setminus S$.

\begin{center}
\psfrag{S1m}{$\Sigma_1$ mixte}
\psfrag{S2}{$\Sigma_2$}
\psfrag{S}{$\Sigma$}
\psfrag{c}{contractant}
\psfrag{U1}{$U_1$}
\includegraphics[width=5.5cm]{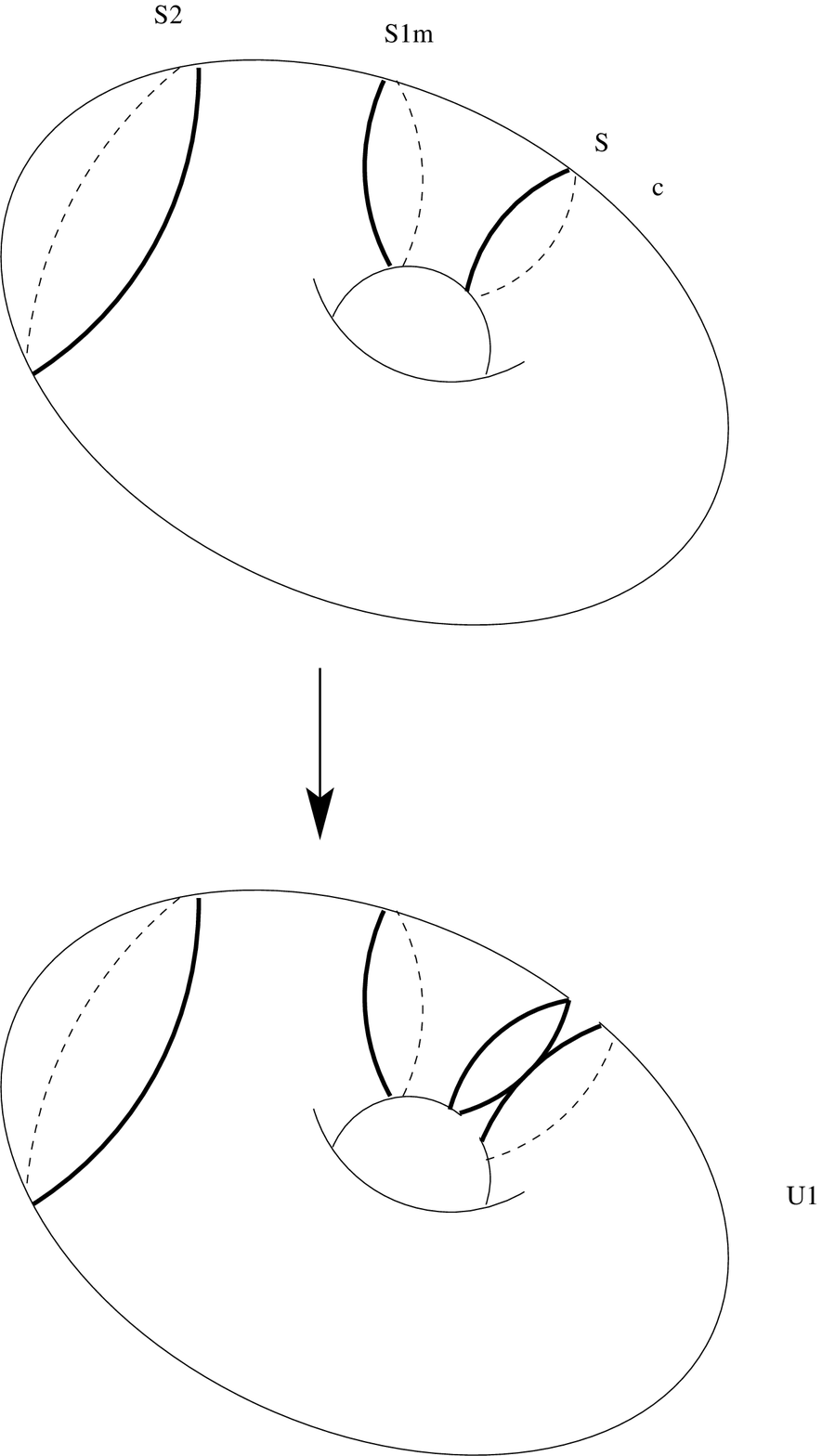}
\end{center}

\bigskip
\begin{center}
Fig.~7: ~ La construction du coeur : 3\`eme possibilit\'e.
\end{center}

{\bf 4\`eme possibilit\'e.} Si les {\bf deux surfaces sont \`a voisinage dilatant}, alors on introduit la vari\'et\'e \`a bord $N^m:=M^m\setminus \Sigma_1\cup\Sigma_2$. Comme $M^m$ est connexe il existe une composante connexe $N^n_1$ de $N^n$ dont le bord contient \`a la fois une copie de $\Sigma_1$ et une copie de $\Sigma_2$. Dans cette composante on minimise l'aire dans la classe d'homologie de $\Sigma_1$. Cette classe d'homologie ne peut pas \^etre nulle dans $H_{m-1}(N^m;{\Z})$ car cette classe est donn\'ee par une union de composantes connexes de $\p N^m_1$ non \'egale \`a la totalit\'e du bord $\p N^n_1$. Le th\'eor\`eme~\ref{reg} combin\'e avec le principe du maximum fort de White (\cite{Whi2} page 420) sur les varifolds d'aire minimale donne l'existence d'une hypersurface $\Sigma$ de $N^m_1$ minimale plong\'ee $S$ et ferm\'ee r\'ealisant cette classe. Le r\'esultat de White donne en particulier que chaque composante connexe de $\Sigma$ est soit totalement incluse dans $\p N^m_1$ soit disjointe
de $\p N^m_1$. Comme ni $\Sigma_1$ ni $\Sigma_2$ n'est localement un minimum de l'aire dans $N^n$, la surface $\Sigma$ ne peut \^etre compl\`etement incluse dans $\p N^m_1$. Il existe donc une composante connexe $\ti{\Sigma}$ de $\Sigma$ de voisinage contractant strictement incluse dans $N^m_1$. On d\'efinit alors $U_1$ comme \'etant une composante connexe de $N^n\setminus \ti{\Sigma}$.

\begin{center}
\psfrag{S1}{$\Sigma_1$}
\psfrag{S2}{$\Sigma_2$}
\psfrag{wS}{$\widetilde{\Sigma}$}
\psfrag{c}{contractant}
\psfrag{U}{$U$}
\includegraphics[width=8cm]{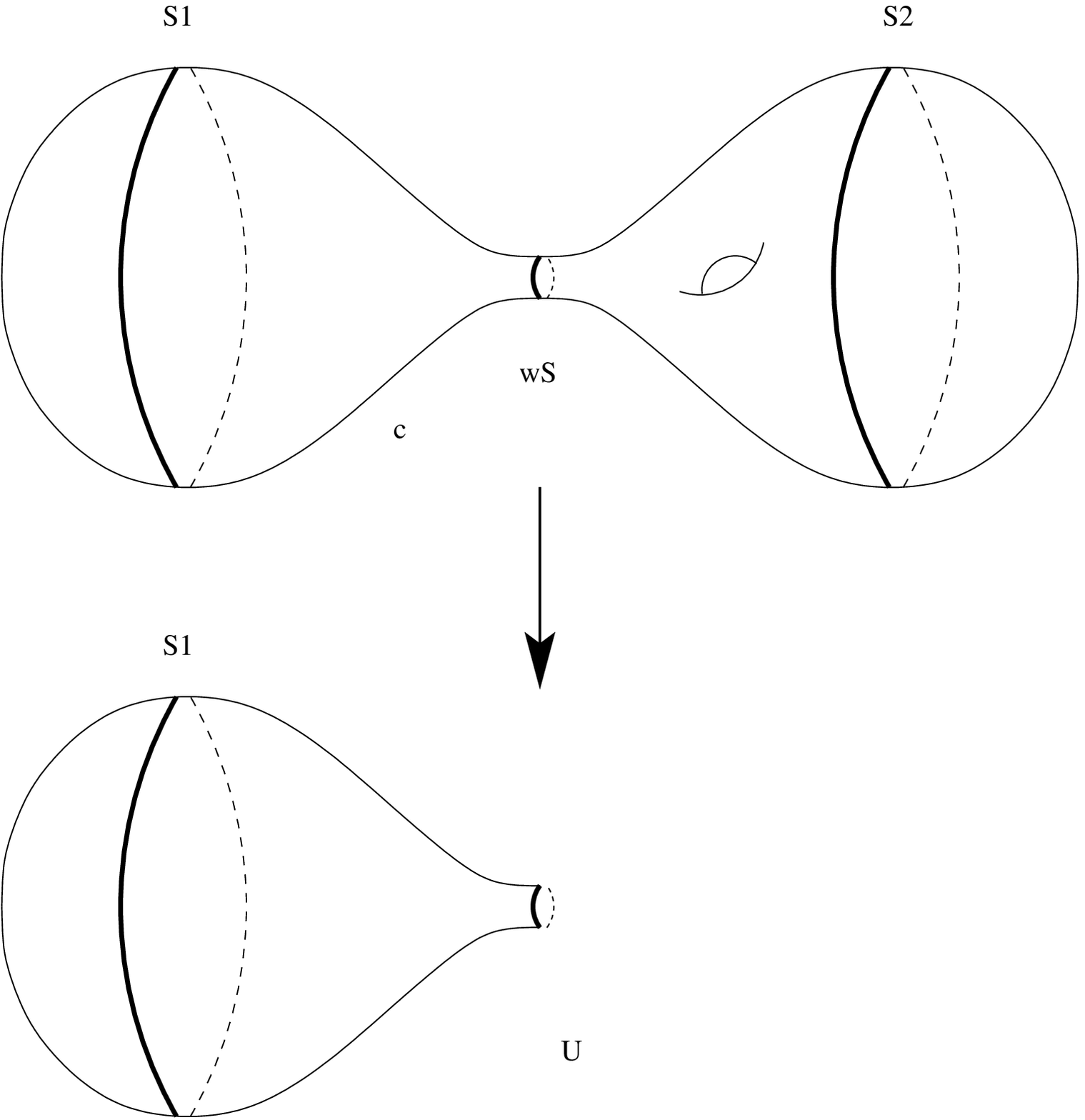}
\end{center}

\bigskip
\begin{center}
Fig.~8 : ~ La construction du coeur : 4\`eme possibilit\'e.
\end{center}

Rempla\c cons maintenant $M^m$ par un ouvert arbitraire $U$ pour lequel  chaque composante du bord est une surface minimale plong\'ee ferm\'ee \`a voisinage contractant.
Si $U$ ne v\'erifie pas la propri\'et\'e de Frankel alors dans chacun des quatre cas ci-dessus on produit un nouvel ouvert $\ti{U}$ strictement inclus dans $U$ et dont chaque composante est une surface minimale \`a bord de voisinage contractant
La proc\'edure  s'it\`ere sans probl\`eme et on produit une suite d\'ecroissante $U_k$ de sous ouverts de $M^m$  dont chaque composante du bord est une surface minimale plong\'ee ferm\'ee  de voisinage contractant. Chacune des \'etape ``consomme'' une nouvelle surface minimale plong\'ee et ferm\'ee de $M^m$, et comme nous supposons qu'il n'en existe qu'un nombre fini nous avons \'etabli le r\'esulat suivant.

\begin{prop}
\label{sous-frankel}
Soit $(M^m,g)$ une vari\'et\'e riemannienne ferm\'ee quelconque de dimension $3\le m\le 7$, si $(M^m,g)$ ne poss\`ede qu'un nombre fini d'hypersurfaces minimales plong\'ees ferm\'ees et g\'eom\'etriquement distinctes alors il existe un ouvert $U$ non vide de $M^m$ tel que les composantes connexes de $\p U$ sont toutes des surfaces minimales plong\'ees ferm\'ees de $M^m$ et \`a voisinage contractant dans $U$.
Toutes les surfaces minimales plong\'ees et ferm\'ees dans l'ouvert $U$ ont un voisinage dilatant.
Par ailleurs, $U$ v\'erifie la propri\'et\'e de Frankel. Un tel ouvert s'appelle un {\bf coeur} de $M^m$.
\end{prop}

On d\'emontre sans trop d'efforts, en  utilisant des arguments semblables \`a ceux utilis\'es pour \'etablir la proposition pr\'ec\'edente, le lemme suivant.
\begin{lemm}
\label{lm-aire}
Soit $U$ un coeur de $M^m$ et $(\Sigma_i)_{i=1\cdots n}$ les composantes connexes du bord de $U$. Alors pour toute hypersurface minimale $\Sigma$ \`a deux faces de $U$ on a
\be
\label{aire1}
\mbox{Vol}(\Sigma)>\max_{i=1,\cdots, n}\mbox{Vol}(\Sigma_i)\quad.
\ee
Pour toute hypersurface \`a une face on a
\be
\label{aire2}
2\,\mbox{Vol}(\Sigma)>\max_{i=1,\cdots, n}\mbox{Vol}(\Sigma_i)\quad.
\ee
\end{lemm}

Afin de compl\'eter la vari\'et\'e \`a bord  r\'ealis\'ee par un coeur donn\'e $U$ on adjoint des cylindres infinis $\Sigma_i\times [0,+\infty)$ \'equip\'es de la m\'etrique produit $g_i+dt^2$ o\`u $\Sigma_i$ sont les composantes connexes du bord de $U$ et $g_i$ la m\'etrique induite par $g$ sur chacune de ces composantes (fig. 9). La vari\'et\'e riemannienne compl\`ete \`a bouts cylindriques obtenue est not\'ee $({\mathcal U},h)$.
Il est important de remarquer \`a ce stade que la m\'etrique $h$ n'est a-priori que Lipschitzienne.

\begin{center}
\psfrag{S1R}{\textcolor{red}{$\Sigma_1 \times {\R}_+$}}
\psfrag{S2R}{\textcolor{red}{$\Sigma_2 \times {\R}_+$}}
\psfrag{US}{$\partial U = \bigsqcup \Sigma_i$}
\psfrag{S1}{$\Sigma_1$}
\psfrag{S2}{$\Sigma_2$}
\psfrag{cU}{le coeur $U$}
\includegraphics[width=10cm]{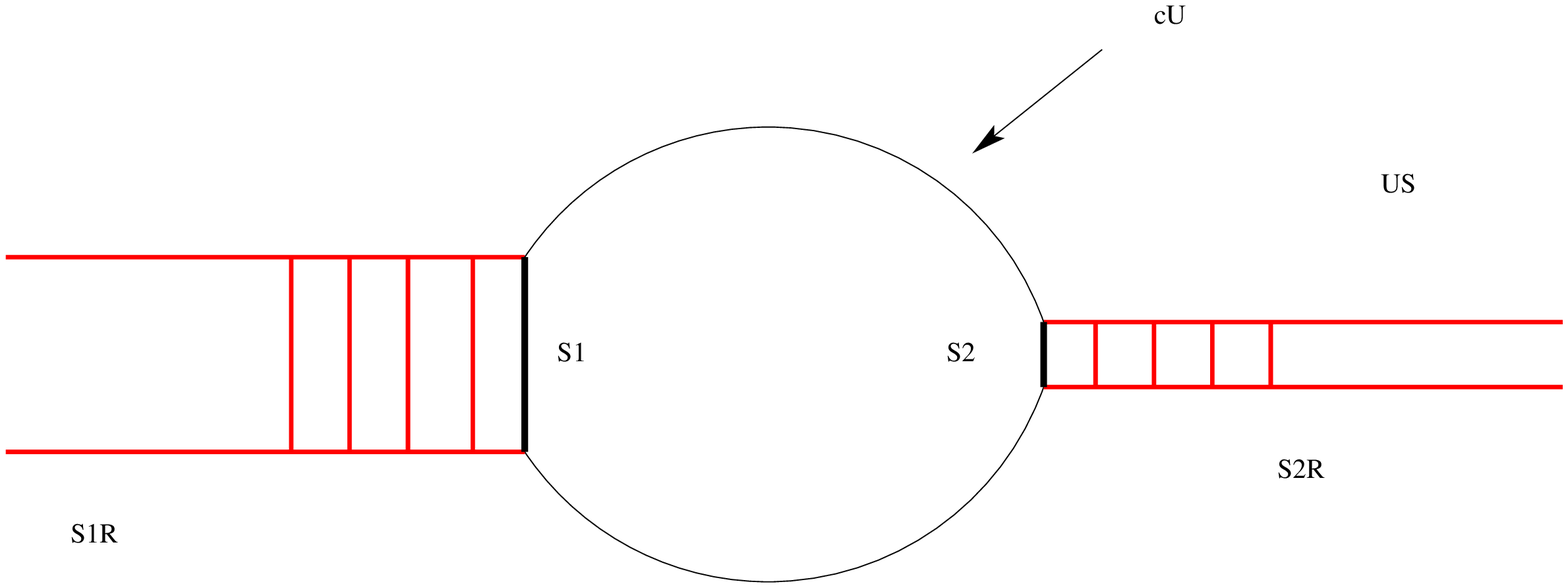}
\end{center}

\begin{center}
${\mathcal U} =$ le coeur $U \;\cup$ \textcolor{red}{les cylindres $\Sigma_i \times {\R}_+$}
\end{center}

\bigskip
\begin{center}
Fig.~9: ~ La compl\'etude du coeur par des cylindres.
\end{center}

Bien que $({\mathcal U},h)$ soit non compacte on peut n\'eanmois d\'efinir son {\bf spectre volumique} $(\la_k(\mathcal U))_{k\in {\N}}$. Pour d\'efinir $\la_k(\mathcal U)$ on observe que pour toute paire $B\subset A$
de sous-ouverts \`a fermeture compacte et \`a bord r\'egulier on a\footnote{La restriction \`a $B$ d'un balayage de $A$ est un balayage de $B$ (modulo le fait que l'op\'eration de restriction n'envoie pas n\'ecessairement  $Z_{m-1,rel}(A,\p A;{\Z}_2)$ sur $Z_{m-1,rel}(B,\p B;{\Z}_2)$ mais seulement $Z_{m-1}(A,\p A;{\Z}_2)$ sur $Z_{m-1}(B,\p B;{\Z}_2)$. Cette difficult\'e peut \^etre r\'esolue \`a l'aide du lemme 2.15 de \cite{LMN} \`a condition que le bord de $B$ soit suffisamment r\'egulier. Nous ignorons ce point technique dans les lignes ci-dessous). Si $j$ est  l'op\'erateur canonique de restriction  des chaines de  $A$ \`a $B$ , 
$j \ : Z_{m-1,rel}(A,\p A;{\Z}_2)\rightarrow Z_{m-1,rel}(B,\p B;{\Z}_2)$, on a donc 
$j^\ast \lambda_B= \lambda_A$. Si $\Psi$ est un $k-$balayage de $A$ on a donc $0\ne\Psi^\ast\lambda_A^k=\Psi^\ast j^\ast \lambda_B^k$ par ailleurs
\[
\max_{x\in X}{\mathbf M}(j\circ\Psi(x))=\max_{x\in X}{\mathbf M}(\Psi(x)\cap B)\le  \max_{x\in X}{\mathbf M}(\Psi(x))
\]
ce qui donne (\ref{inclus}).}
\be
\label{inclus}
\la_k(A)\ge \la_k(B)
\ee
Pour toute suite croissante d'ouverts \`a fermeture compacte et bords r\'eguliers $K_1\subset K_2\cdots$ telle que $\cup_{i\in {\N}} K_i={\mathcal U}$ on d\'efinit
\[
\la_k({\mathcal U}):=\lim_{i\rightarrow+\infty}\la_k(K_i)\quad.
\]
On v\'erifie sans trop d'efforts que cette limite est finie et est ind\'ependante du choix de la suite $K_i$. La prochaine \'etape de l'argument consiste \`a \'etablir le r\'esultat suivant.
\begin{lemm}
\label{lm-width-coeur} Soit $U$ un coeur d'une vari\'et\'e riemannienne $(M^m,g)$ et soit ${\mathcal U}$ la vari\'et\'e riemannienne compl\`ete obtenue \`a partir de $U$ en ajoutant des bouts cylindriques
$\p U\times {\R}_+^\ast$ \'equip\'es de la m\'etrique produit. 

Soit $\Sigma_1$ la composante connexe de $\p U$ de volume maximal. On a d'une part pour tout $k=1,2,3\cdots$
\be
\label{croiss-spec}
\la_{k}({\mathcal U})+\mbox{Vol}(\Sigma_1)\le \la_{k+1}({\mathcal U})\quad.
\ee
Par ailleurs, il existe $C({\mathcal U})>0$ telle que pour tout $k=1,2,3\cdots$
\be
\label{bornes-spec}
k\,\mbox{Vol}(\Sigma_1)\le \la_{k}({\mathcal U})\le k\,\mbox{Vol}(\Sigma_1)+C({\mathcal U})\ k^{\frac{1}{m+1}}\quad.
\ee
\hfill $\Box$
\end{lemm}
\begin{rema}
\label{rm-width-coeur}
Il est int\'eressant d'observer que la croissance lin\'eaire du spectre volumique du coeur, \'etendu par des bouts cylindriques infinis, et donn\'ee par (\ref{bornes-spec}) contraste avec la croissance sous-lin\'eaire (\ref{5}) du coeur lui m\^eme ou de la vari\'et\'e ferm\'ee $(M^m,g)$ dont il est issu. \hfill $\Box$
\end{rema}
La preuve du lemme~\ref{bornes-spec} peut se comprendre ainsi. Tout d'abord, il n'est pas difficile de se convaincre intuitivement (voir les arguments complets dans \cite{Son}) que la premi\`ere largeur de Gromov, c'est \`a dire l'aire maximale d'un balayage, d'un cylindre suffisament long est donn\'ee par sa section. On a donc l'existence de $R_0$ tel que pour tout $R\ge R_0$
\[
\la_1(\Sigma_1\times [0,R], g_1+dt^2)=\mbox{Vol}(\Sigma_1)\quad.
\]
Soit $p$ un point de ${\mathcal U}$. Pour tout $\ep$ et $R$ suffisamment grand $\la_k(B_R(p))\ge \la_k({\mathcal U})-\ep$, o\`u $B_R(p)$ est la boule g\'eod\'esique de centre $p$ et de rayon $R$. On choisit $\rho$ assez grand tel que $B_R(p)\cap \Sigma_1\times[\rho,\rho+R]=\emptyset$. Les propri\'et\'es fondamentales du cup produit donnent de fa\c con analogue \`a la d\'emonstration de (\ref{6}) (fig.10)
\[
\begin{array}{c}
\ds\la_{k+1}({\mathcal U})\ge \la_{k+1}\lf(B_R(p)\cup\Sigma_1\times[\rho,\rho+R]\rg)\\[5mm]
\ds\quad\quad\quad\quad\quad\quad\quad\quad\ge \la_k(B_R(p))+\la_1(\Sigma_1\times [\rho,\rho+R])\ge \la_k({\mathcal U})-\ep+\mbox{Vol}(\Sigma_1)\quad.
\end{array}
\]
On en d\'eduit (\ref{croiss-spec}). L'in\'egalit\'e (\ref{inclus}) et la d\'efinition du spectre volumique de ${\mathcal U}$ donnent  par ailleurs
\[
\la_1({\mathcal U})\ge \la_1(\Sigma_1\times [0,R_0], g_1+dt^2)=\mbox{Vol}(\Sigma_1)
\]
En combinant cette derni\`ere in\'egalit\'e avec (\ref{croiss-spec}) on obtient la minoration annonc\'ee de $\la_k({\mathcal U})$ dans (\ref{bornes-spec}).

\begin{center}
\psfrag{S1R}{\textcolor{red}{$\Sigma_1 \times {\R}_+$}}
\psfrag{S2R}{\textcolor{red}{$\Sigma_2 \times {\R}_+$}}
\psfrag{US}{$\partial U = \bigsqcup \Sigma_i$}
\psfrag{S1}{$\Sigma_1$}
\psfrag{S2}{$\Sigma_2$}
\psfrag{cU}{le coeur $U$}
\psfrag{BR}{\textcolor{green}{$B_R(p)$}}
\psfrag{S1x}{\textcolor{green}{$\Sigma_1 \times [p; p + R_0]$}}
\psfrag{xp}{\textcolor{green}{$x_p$}}
\includegraphics[width=13cm]{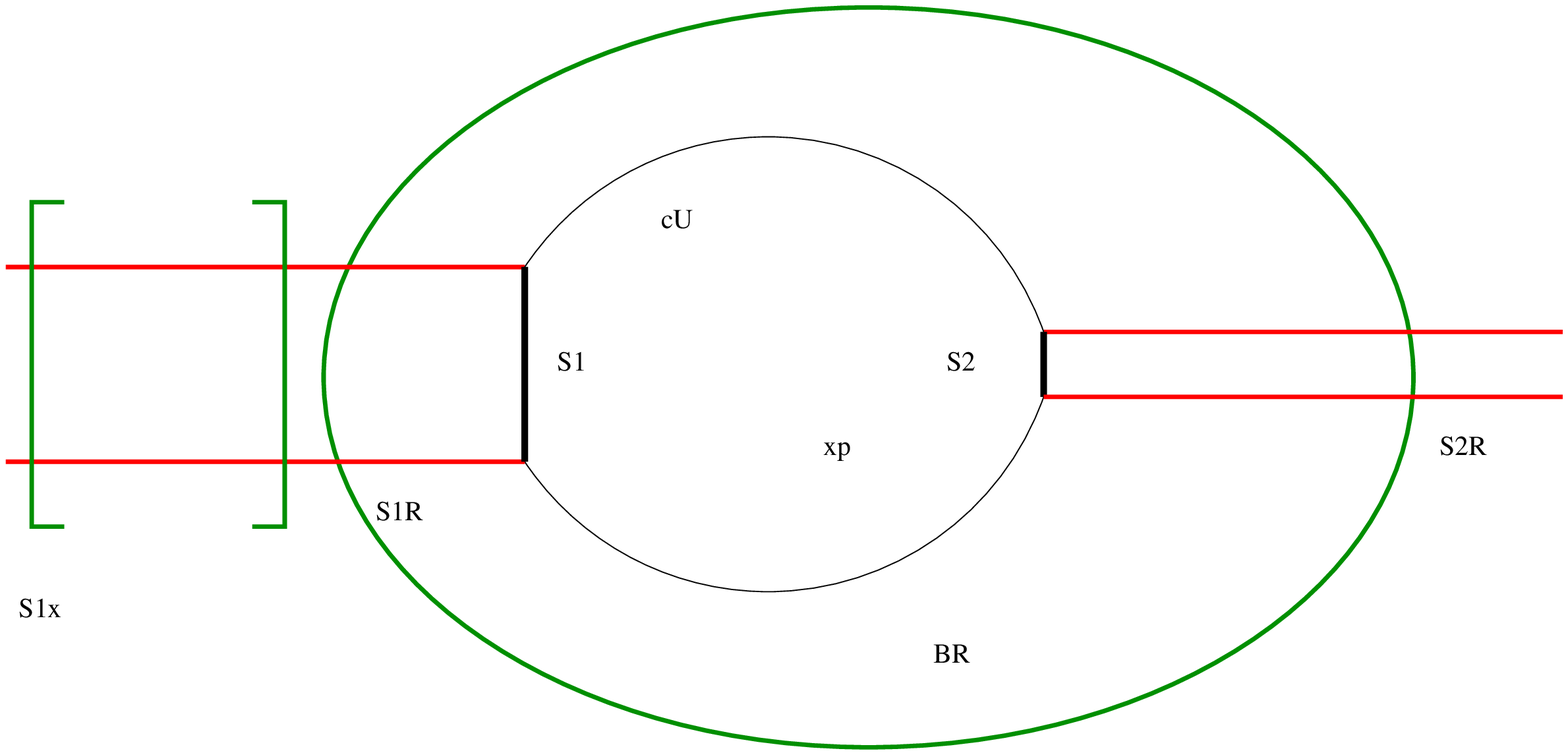}
\end{center}

\begin{center}
$\lambda_1\big(\Sigma_1 \times [p; p + R_0]\big) + \lambda_k \big(B_R (p)\big) \le \lambda_{k+1}({\mathcal U})$
\end{center}

\bigskip
\begin{center}
Fig.~10: ~ La minoration du spectre volumique
\end{center}

La majoration est obtenue en ``recollant'' un $k-$balayage de l'union des cylindres $\p U\times {\R}_+^\ast$ qui contribue \`a $k\times \mbox{Vol}(\Sigma_1)$ avec un
$k-$balayage du coeur donn\'e par (\ref{5}) et dont le co\^ut est contr\^ol\'e par $C\ k^{1/(m+1)}$. Ce recollement de $k-$balayage n'est pas ais\'e et requiert une technique
introduite dans \cite{LMN}. On obtient la majoration annonc\'ee de $\la_k({\mathcal U})$ dans (\ref{bornes-spec}).

On aborde maintenant la partie la plus d\'elicate de la preuve. Il s'agit de d\'emontrer que chaque valeur du spectre de ${\mathcal U}$ est r\'ealis\'ee par une surface minimale plong\'ee
\underbar{incluse dans l'int\'erieur de $U$} et \'equip\'ee d'une multiplicit\'e enti\`ere. Pr\'ecisemment Antoine Song \'etablit le r\'esultat suivant.
\begin{lemm}
\label{lm-real}
Soit $(\lambda_k({\mathcal U}))_{k\in {\N}^\ast}$ le spectre volumique de la compl\'etude cylindrique d'un coeur $U$ d'une vari\'et\'e riemannienne $(M^m,g)$, alors, pour tout $k=1,2,3\cdots$ il existe
une surface minimale plong\'ee ferm\'ee $S_k$ de $U$, incluse dans $U\setminus \p U$, et une multiplicit\'e $m_k\in {\N}^\ast$ tels que
\be
\label{reali}
\la_k({\mathcal U})=m_k\, \mbox{Vol}(S_k)\quad.
\ee
\hfill $\Box$
\end{lemm}
Le sh\'ema de la preuve est assez simple. Pour une suite d'ouverts  \`a fermeture compacte bien choisis de ${\mathcal U}$ et convergeants vers ${\mathcal U}$ on effectue une op\'eration de minmax sur chacun d'entre eux afin de mettre en \'evidence des surfaces minimales r\'ealisants leurs spectres volumiques. Puis on passe \`a la limite sur ces surfaces minimales.

 Avant m\^eme de se poser la question du passage \`a la limite et pourquoi
les surfaces limites sont n\'ecessairement incluses dans l'int\'erieur de $U$, tel que l'annonce le th\'eor\`eme, on doit constater que les op\'erations de minmax envisag\'ees sortent du cadre habituel du th\'eor\`eme~\ref{almpitt} d'Almgren et Pitts et cela pour deux raisons : 1) les ouverts d\'efinissants 
le spectre ont des bords 2) la m\'etrique $h$ \'equipant ${\mathcal U}$ n'est que Lipschitz sur ces ouverts et donc pas suffisament r\'eguli\`ere pour pouvoir se reposer sur la th\'eorie classique.

On modifie alors la m\'etrique $h$ de ${\mathcal U}$ dans un voisinage contenant strictement $\p U\times {\R}_+^\ast$ de la fa\c con suivante. Comme chaque composante de $\p U$ est une surface minimale contractante
on peut modifier $h$ en une m\'etrique r\'eguli\`ere $h_\ep$, convergeant en norme $C^l_{loc}({\mathcal U}\setminus \p U)$ pour tout $l\in {\N}$,  et tel qu'un voisinage de chaque cylindre $\p \Sigma_i\times {\R}_+^\ast$ poss\`ede un feuilletage r\'egulier dont le champs de vecteur courbure moyenne pointe vers l'infini du cylindre (un tel feuilletage est appel\'e {\bf convexe en moyenne})\footnote{En particulier les surfaces $\Sigma_i$ ne peuvent plus \^etre minimale pour $h_\ep$ \`a cause du principe du maximum.}. On approche par ailleurs
$({\mathcal U},h_\ep)$ par une suite d'ouverts $({\mathcal U}_\ep,h_\ep)$ \`a fermeture compacte et de bords r\'eguliers en tronquant simplement les cylindres le longs de feuilles de plus en plus \'eloign\'ees du coeur (figure 11).

\medskip

\begin{center}
\psfrag{S1R+}{\textcolor{red}{$\simeq \Sigma_1 \times [0,\ep^{-1})$}}
\psfrag{S2R+}{\textcolor{red}{$\simeq \Sigma_2 \times [0,\ep^{-1})$}}
\psfrag{lc}{le coeur}
\psfrag{U}{$U$}
\psfrag{Uh}{$({\mathcal U}_\varepsilon, h_\varepsilon)$}
\psfrag{feuilletage}{\textcolor{red}{Le feuilletage convexe en moyenne des voisinages des cylindres}}
\includegraphics[width=10cm]{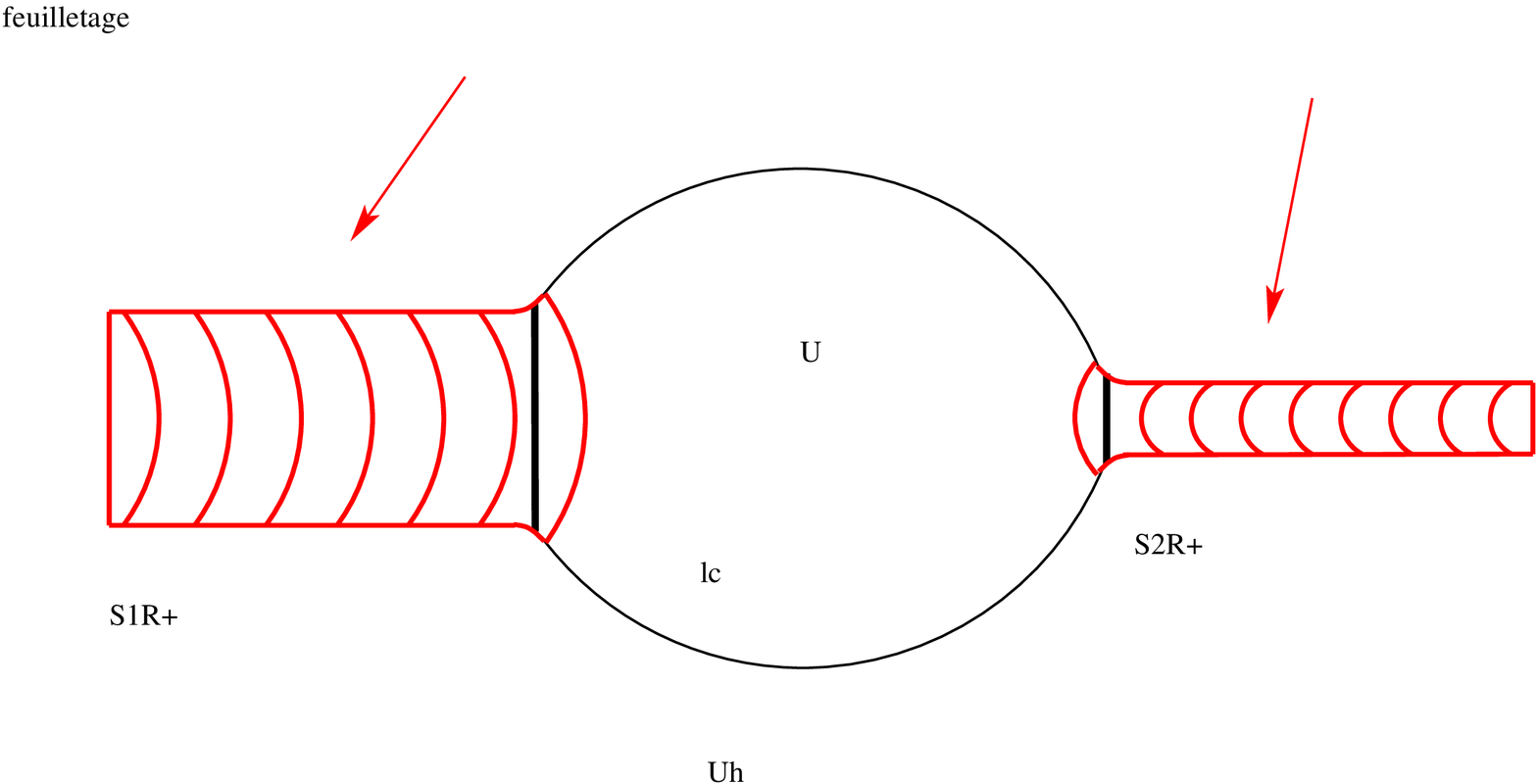}
\end{center}

\bigskip
\begin{center}
Fig.~11: ~ La m\'etrique approxim\'ee de la compl\'etude du coeur.
\end{center}

Un r\'esultat - non encore publi\'e - de Martin Li et Xin Zhou \cite{LZ} permet d'\'etendre la th\'eorie d'Almgren-Pitts-Marques-Neves au cas des vari\'et\'es compactes \`a bord. D'apr\`es ce travail technique et difficile il existe donc pour tout $k$ et tout $\epsilon>0$ une famille finie de surfaces $\{S_1^{\,\ep,k}\cdots S_{n_{\ep,k}}^{\,\ep,k}\}$ compactes plong\'ees de $({\mathcal U}_\ep,h_\ep)$ qui sont minimales dans l'int\'erieur de ${\mathcal U}_\ep$ et \`a {\bf bord libre}\footnote{Les surfaces minimales plong\'ees \`a bord libre sont des surfaces minimales compactes $\Sigma$ de $({\mathcal U}_\ep,h_\ep)$, telles que $\p\Sigma\subset \p {\mathcal U}_\ep$, qui sont continuement diff\'erentiables jusqu'au bord et dont le vecteur normal au plan limite tangent au bord est contenu dans $\p {\mathcal U}_\ep$.} ainsi qu'une famille d'entiers non nuls $m_i^{\ep,k}$ tels que
\[
\la_k({\mathcal U}_\ep,h_\ep)=\sum_{i=1}^{n_{\ep}(k)}m_i^{\ep}(k)\ \mbox{Vol}(S_i^{\,\ep}(k))\quad.
\]
Nous d\'emontrons maintenant que pour tout entier $k$ et tout nombre strictement positif $\ep$ suffisamment petit, toutes les surfaces $S_i^{\,\ep}(k)$ sont incluses  dans un compact, d\'ependant de $k$ mais ind\'ependant de $\epsilon$. Supposons 
dans un premier temps qu'une surface $S_i^{\,\ep}(k)$ poss\`ede une composante connexe $\ti{S}$ strictement incluse dans un cylindre. Il existe donc une feuille du feuilletage la plus proche du coeur et tangente \`a
$\ti{S}$. Cette feuille a le vecteur courbure moyenne pointant vers l'infini du cylindre et donc vers le c\^ot\'e o\`u $\ti{S}$ se trouve. Une application directe du principe du maximum apporte une contradiction.
On en d\'eduit que toute composante connexe $\ti{S}$ de $S_i^{\,\ep}(k)$ doit intersecter le coeur (fig. 12).

\begin{center}
\psfrag{Sk}{$S^\varepsilon_i(k)$}
\psfrag{lc}{le coeur}
\psfrag{U}{$U$}
\psfrag{tmin}{\textcolor{red}{$t_{\min}$}}
\psfrag{possible}{\textcolor{green}{possible}}
\psfrag{impossible}{\textcolor{green}{impossible}}
\includegraphics[width=12cm]{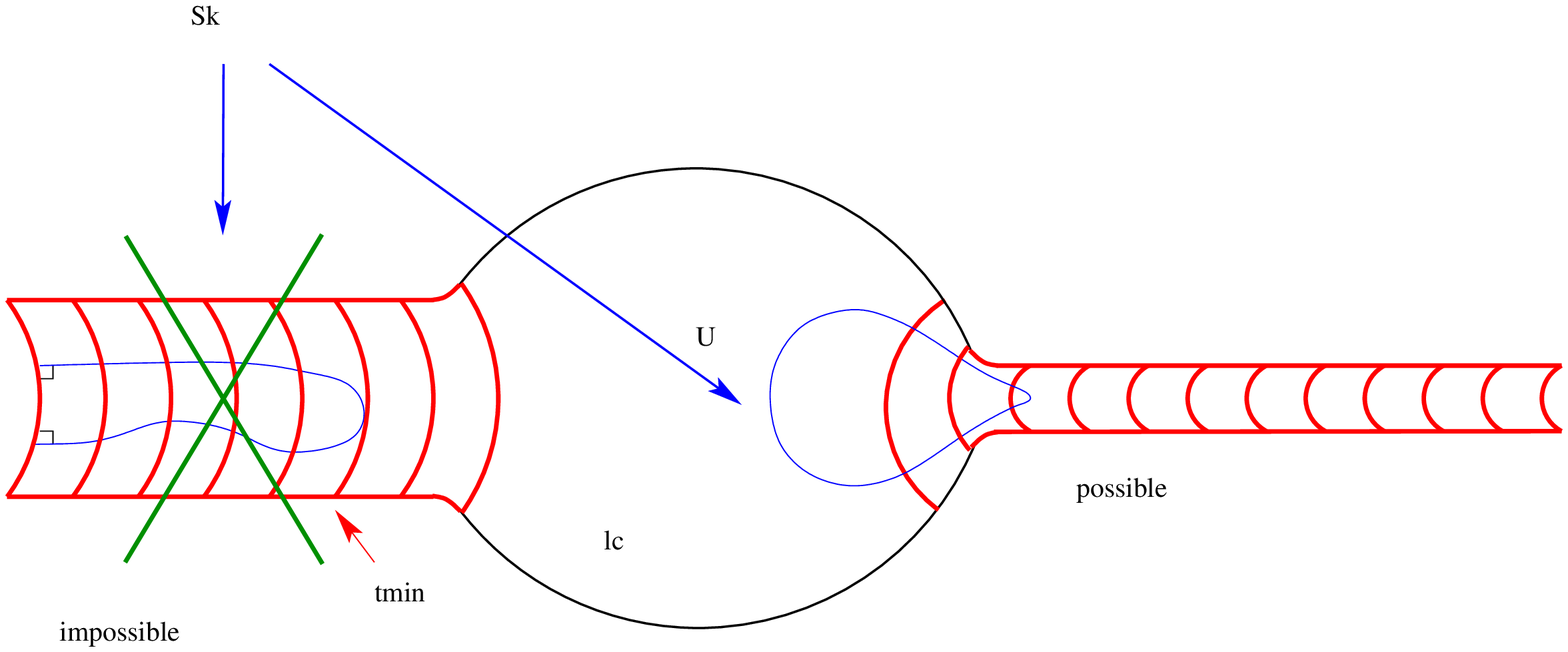}
\end{center}

\bigskip
\begin{center}
Fig.~12: ~ Le feuilletage convexe en moyenne agissant comme ``barri\`ere''
\end{center}

 Si un point d'une telle composante connexe $\ti{S}$ est \'eloign\'e d'une distance $L$ du coeur et poss\`ede un point $p\in \p \Sigma_i\times \{L\}$ ,  de la formule de monotonie sur les surfaces minimales\footnote{La formule de monotonie pour les hypersurfaces minimales dans l'espace espace euclidien ${\R}^m$ par exemple affirme que pour toute surface minimale ferm\'ee $S$ et pour tout point $p$ de ${\R}^m$ la fonction
 \[
 r\ \longrightarrow\ \frac{\mbox{Vol}(S\cap B^m_r(p))}{r^{m-1}}
 \]
est croissante, o\`u $B^m_r(p)$ est la boule euclidienne de centre $p$ et de rayon $r>0$ (\cite{Sim}). Si donc $p$ est un point de $S$, comme la limite en z\'ero de cette fonction est \'egale au volume de $B^{m-1}_1(0)$ c\'est \`a dire \`a $\pi^{(m-1)/2}/\Gamma((m+1)/2)$. On obtient la minoration
\[
\mbox{Vol}(S\cap B^m_r(p))\ge \frac{\pi^{(m-1)/2}}{\Gamma((m+1)/2)}\ r^{m-1}\quad.
\]
Modulo une modification des constantes, une telle minoration est aussi valable pour $r\le 1$ dans une vari\'et\'e riemannienne ferm\'ee donn\'ee ou encore dans une vari\'et\'e riemannienne compl\`ete dont la g\'eom\'etrie est born\'ee comme les cylindres $\Sigma_i\times {\R}$. C'est cette minoration qui est responsable du confinement de $S^\ep_i(k)$ et qui exclut la possibilit\'e pour l'hypersurface $\ti{S}$ de s'\'etendre hors d'un compact dans les cylindres (voir fig. 13). } on d\'eduit que sur chaque segment ${\p \Sigma_i}\times [n, n+1]$ l'aire de $\ti{S}$ est au moins \'egale \`a une constante $C>0$ ne d\'ependant que de $\Sigma_i$.  Donc, on a
$$\mbox{Vol}(\ti{S})\ge C\, L^{m-1}$$ o\`u $C>0$ ne d\'epend que de $U$. Comme $\mbox{Vol}(\ti{S})\le \la_k({\mathcal U}_\ep)$ on peut d\'eduire l'in\'egalit\'e suivante $$L< C^{-1/(m-1)}\ \la^{1/(m-1)}_k({\mathcal U}_\ep)\quad.$$
Donc, les surfaces $S_i^{\,\ep}(k)$ sont incluses  dans un compact, d\'ependant de $k$ mais ind\'ependant de $\epsilon$ (fig. 13).

\begin{center}
\psfrag{Sk}{$S^\varepsilon_i(k)$}
\psfrag{lc}{le coeur}
\psfrag{U}{$U$}
\psfrag{impossible}{\textcolor{green}{impossible}}
\includegraphics[width=12cm]{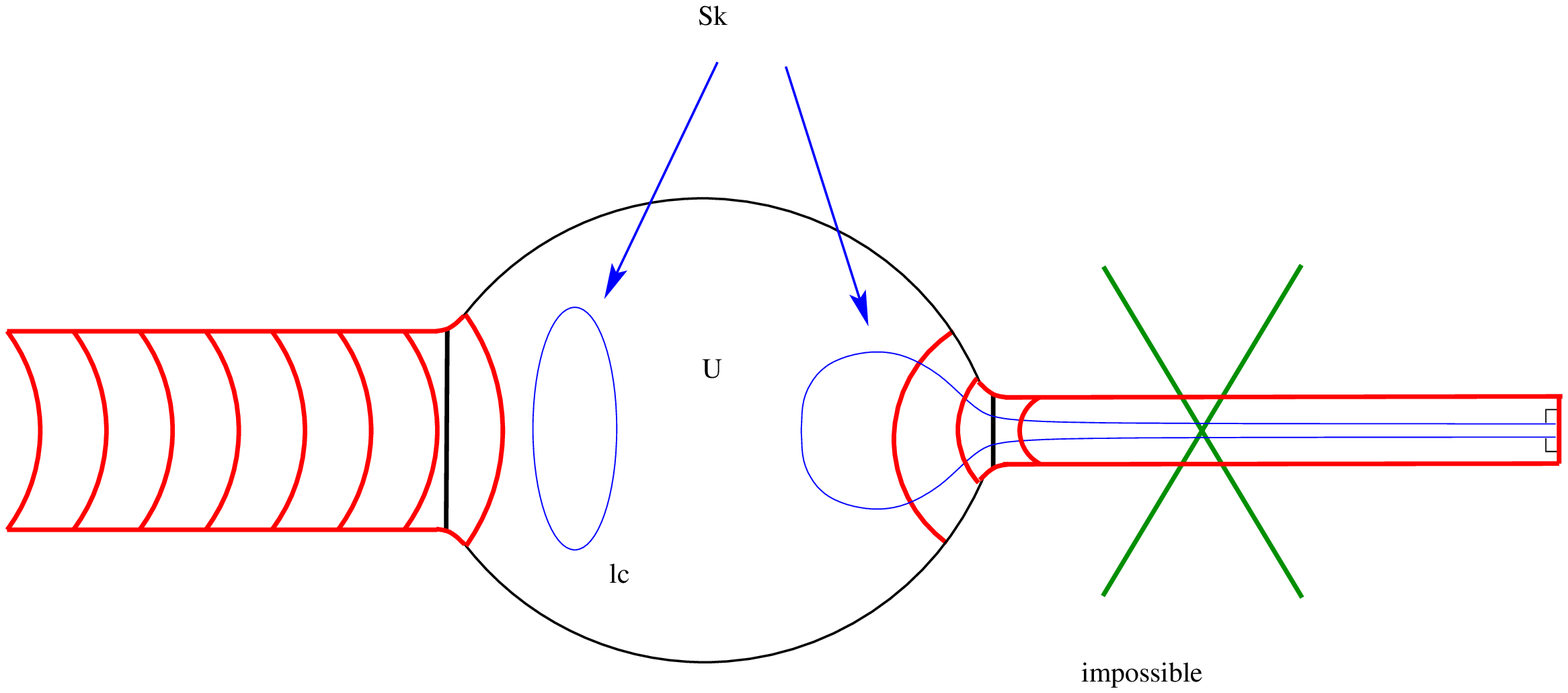}
\end{center}

\bigskip
\begin{center}
Fig.~13: ~ Le confinement des hypersurfaces minimales pr\`es du coeur
\end{center}

Dans \cite{MN1} F.C. Marques et A. Neves d\'emontrent - dans le cas des vari\'et\'es ferm\'ees - que le spectre volumique est effectivement r\'ealis\'e par  une famille de surfaces minimales $S_i^{\,\ep}(k)$ satisfaisant
\be
\label{indice}
\sum_{i=1}^{n_{\ep}(k)}\mbox{Ind}(S_i^{\,\ep}(k))\le k\quad.
\ee
On peut alors invoquer le r\'esultat de compacit\'e de Ben Sharp \cite{Sha} et passer \`a la limite dans les familles de varifolds ${\mathbf v}^{\ep}(k)$ donn\'ees par les sommes des mesures de Radons port\'ees par les surfaces $S_i^{\,\ep}(k)$ 
pour les multiplicit\'es $m_i^{\ep}(k)$. La limite est le varifold ${\mathbf v}_k$ donn\'e par une hypersurface $S_k$ minimale plong\'ee r\'eguli\`ere  en dehors du bord\footnote{Le r\'esultat de Sharp ne peut s'appliquer qu'en dehors de $\p U$ l\`a o\`u la m\'etrique $h_\ep$ converge fortement vers $h$ dans toutes les normes $C^l$ pour $l\in {\N}$} $\p U$ et affect\'ee de multiplicit\'es enti\`ere r\'eguli\`ere en dehors de $\p U$. Cette surface  $S_k$ ne peut inters\'ecter l'int\'erieur des cylindres $\p \Sigma_i\times {\R}^\ast_+$. En effet ceux-ci 
sont feuillet\'es par les surfaces minimales $\p \Sigma_i\times\{t \}$ et en choisissant, pour un cylindre donn\'e, le $t=t_{max}$ maximal d'intersection avec $S_k$, le principe du maximum de White \cite{Whi2} impose qu'une composante 
connexe de $S_k$ doit \^etre \'egale \`a $\p \Sigma_i\times\{t_{max} \}$. Cette derni\`ere affirmation contredit le fait qu'\`a $\ep>0$ fix\'e, toute composante connexe des surfaces $S_i^{\,\ep}(k)$ doive intersecter le coeur.
En effet, la formule de monotonie pr\'ec\'edemment utilis\'ee impose qu'\`a $\ep>0$ toute portion du cylindre entre $t=0$ et $t=t_{max}$ contient une quantit\'e de volume au moins proportionnel \`a la longueur de cette portion, ind\'ependamment de $\ep$. \`A cause de la nature de la convergence\footnote{Il s'agit d'une convergence pour les varifolds et donc en mesure de Radon pour laquelle la masse passe \`a la limite.}, ce volume n'a pu s'\'evanouir dans la section $\p \Sigma_i\times(0,t_{max})$.

On a donc $S_k\subset \ov{U}$ et cette hypersurface r\'ea lise une hypersurface minimale plong\'ee dans sa partie inters\'ectant l'int\'erieur de $U$.

\begin{center}
\psfrag{S1R}{\textcolor{red}{$\Sigma_1 \times {\R}_+$}}
\psfrag{S2R}{\textcolor{red}{$\Sigma_2 \times {\R}_+$}}
\psfrag{Sk}{\textcolor{blue}{$S_k$}}
\psfrag{cU}{le coeur $U$}
\psfrag{tmax}{\textcolor{red}{$t_{\max}$}}
\psfrag{impossible}{\textcolor{green}{impossible}}
\psfrag{zone}{\textcolor{blue}{zone de contact $\partial U \cap S_k$}}
\includegraphics[width=12cm]{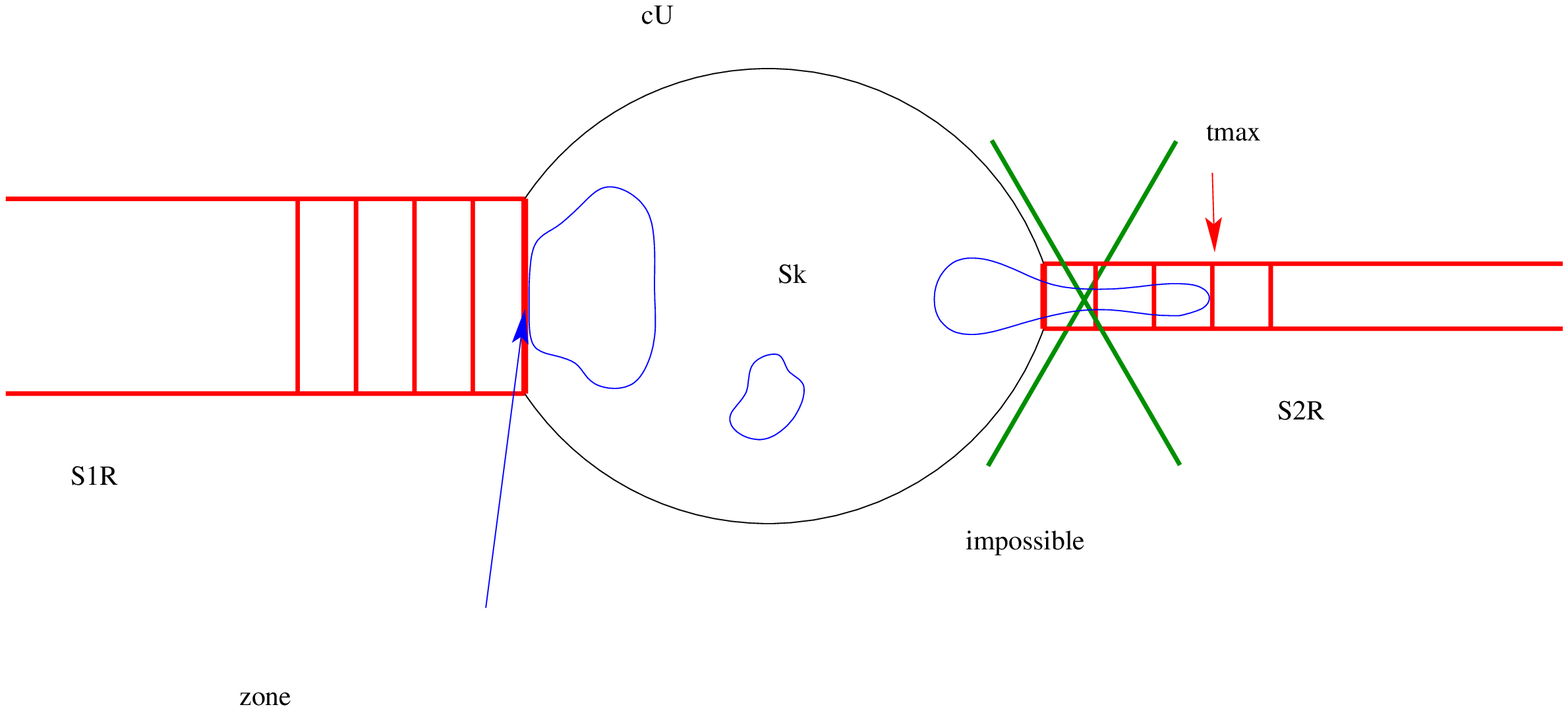}
\end{center}

\bigskip
\begin{center}
Fig.~14: ~ Le passage \`a la limite $\varepsilon \rightarrow 0$.
\end{center}

 Le probl\`eme vient de la possible existence d'une ``zone de contact'' $S_k\cap \p{U}$ (fig.14). En effet, le long de $\p U$  la convergence de la m\'etrique $h_\ep$ vers $h$ n'est pas assez ``forte'' afin de conclure que $S_k$ est aussi minimale et plong\'ee  dans un voisinage de $\p U$. Ce fait est \'etabli au moyen d'un calcul tr\`es d\'elicat\footnote{C'est peut-\^etre l'\'etape de l'article de Song la plus technique et difficile \`a v\'erifier pour un lecteur exigeant.} qui d\'emontre que le varifold
limite ${\mathbf v}_k$ est bien stationnaire dans $(U,g)$.

En invoquant \`a nouveau le principe du maximum de White on en d\'eduit que ou bien $S_k\cap \p{U}$ est donn\'ee par une union de composantes connexes de $\p U$, les  hypersurfaces $\Sigma_i$, ou bien
$S_k\cap \p{U}=\emptyset$. Le premier cas est exclus pour les m\^emes raisons\footnote{En effet, il faut se souvenir \`a ce stade que pour $\ep>0$ le feuilletage convexe moyen existe aussi dans un voisinage de $\p U$ dans $U$ - et pas seulement dans les cylindres - et que donc toute composante connexe des surfaces $S_i^{\,\ep}(k)$ doit  inters\'ecter un point du coeur \underbar{en dehors} d'un voisinage de $\p U$. C'est un passage important de la preuve qui n\'ecessite que toutes les surfaces $\p \Sigma_i$ soient contractantes et qui exclut en particulier la possibilit\'e que le varifold limite  ${\mathbf v}_k$ soit \'egal \`a $k-$fois la mesure de Radon donn\'ee par $\Sigma_1$.}
 qui nous ont permis d'exclure que $S_k$ ne peut contenir une ``tranche de cylindre'' $\Sigma_i\times \{t\}$.

La surface limite r\'ealisant $\la_k({\mathcal U},h)$ est donc strictement incluse dans le coeur. Comme ce dernier v\'erifie la propri\'et\'e de Frankel on en d\'eduit le lemme~\ref{lm-real}.

En supposant donc que le coeur ne poss\`ede qu'un nombre fini de surfaces minimales plong\'ees, sachant qu'elles sont toutes de volume strictement sup\'erieur\footnote{Au cours de la preuve on d\'emontre que les surfaces  \`a une face sont toutes r\'ealis\'ees un nombre pair de fois.} \`a $\mbox{Vol}(\Sigma_1)$ (lemme~\ref{lm-aire}) et sachant que $\la_{k+1}({\mathcal U},h)\ge \la_{k}({\mathcal U},h)+k\,\mbox{Vol}(\Sigma_1)$ un lemme de combinatoire
donne alors l'existence d'un nombre positif $\al>0$ tel que pour $k$ suffisamment grand
\[
\la_{k}({\mathcal U},h)\ge (1+\al)\ k\,\mbox{Vol}(\Sigma_1)
\]
Ceci contredit la minoration donn\'ee par (\ref{bornes-spec}). On en conclut que $M^m$ poss\`ede un nombre infini d'hypersurfaces minimales plong\'ees ferm\'ees g\'eom\'etriquement distinctes.
Le th\'eor\`eme~\ref{th-song} est d\'emontr\'e.

\section{Conclusion}

Nous avons d\^u faire des choix pour cette pr\'esentation. Certaines questions importantes autour de la probl\'ematique de la r\'ealisation du spectre volumique par des surfaces minimales ont germ\'e
\`a partir de la th\'eorie d'Almgren-Pitts. Une des question est celle de la valeur de l'indice de Morse de chacune de ces surfaces. Nous avons mentionn\'e bri\`evement un r\'esultat partiel de Marques et Neves allant dans cette direction
et ceci \`a l'occasion de passages \`a la limite dans la partie VIII de cet expos\'e. Ce r\'esultat (\ref{indice}) a \'et\'e compl\'et\'e r\'ecemment par deux contributions importantes. La premi\`ere est une pr\'ebublication \`a nouveau de Marques et Neves  \cite{MN2} qui \'etablit des conditions suffisantes afin d'avoir une \'egalit\'e dans (\ref{indice}) lorsque les multiplicit\'es des hypersurfaces minimales r\'ealisant une valeur propre donn\'ee du spectre sont toutes \'egales \`a un. La deuxi\`eme contribution est une pr\'epublication de Xin Zhou \cite{Zho} annon\c cant la preuve, dans le cas des m\'etrique g\'en\'eriques,  d'une conjecture importante, {\bf la conjecture de multiplicit\'e un},  suivant laquelle les multiplicit\'es des hypersurfaces minimales r\'ealisant une valeur propre donn\'ee du spectre sont toutes \'egales \`a un.

Dans une autre direction, une fois la conjecture de Yau d\'emontr\'ee, il est int\'eressant d'explorer la ``r\'epartition'' des surfaces minimales plong\'ees d'une vari\'et\'e donn\'ee. Un r\'esultat tr\`es int\'eressant de Marques, Neves et Song \cite{MNS} affirme que dans le cas de m\'etriques g\'en\'eriques il existe une famille d\'enombrable d'hypersurfaces minimales plong\'ees et ferm\'ees satisfaisant une propri\'et\'e d'{\bf \'equidistribution}.

La th\'eorie d'Almgren-Pitts n'est pas l'approche exclusive jusqu'ici adopt\'ee  pour la construction d'hypersurfaces minimales en vue de r\'esoudre la conjecture de Yau.
La  m\'ethode dite des {\bf ensembles de niveau} a r\'ecemment port\'e ses fruits dans un article remarquable d'Otis Chodosh et de Christos Mantoulidis \cite{CM}. Ils  donnent une 
d\'emonstration tout \`a fait originale du th\'eor\`eme~\ref{infweyl} dans le cas de la dimension 3 en \'etudiant la concentration des ensembles de niveau de la fonctionnelle d'{\bf Allen-Cahn}\footnote{Les \'equations d'Allen-Cahn,  sont des \'equations mod\'elisant des ph\'enom\`enes de transition de phase en physique des mat\'eriaux.}  lorsque le param\`etre de couplage entre l'\'energie \'elastique et le potentiel \`a ``deux puits'' tend vers l'infini. Ce projet d'utiliser la fonctionnelle d'Allen-Cahn pour construire par minmax des hypersurfaces minimales  remonte \`a des travaux de John E. Hutchinson and Yoshihiro Tonegawa \cite{HT} avec des contributions r\'ecentes  en particulier de Marco Guaraco et Pedro Gaspar (\cite{Gua}, \cite{GG}). Une pr\'esentation d\'etaill\'ee de cette approche qui utilise des r\'esultats tr\`es fins d'analyse d'\'equations aux d\'eriv\'ees partielles elliptiques semi-lin\'eaires scalaires nous emmenerais bien au-del\`a du cadre de cet expos\'e.

Finallement, il est naturel de se poser la question de l'extension de la conjecture de Yau en codimension quelconque\footnote{Comme on le mentionnait dans la premi\`ere partie de l'expos\'e, cette question est encore ouverte d\'ej\`a dans le cas tr\`es particulier des g\'eod\'esiques ferm\'ees sur les sph\`eres  arbitraires de dimension sup\'erieures \`a deux.} : 

{\it Question :  Dans une vari\'et\'e riemannienne ferm\'ee  arbitraire existe t'il un nombre infini de sous vari\'et\'es (ou m\^eme d'immersions partiellement r\'eguli\`eres) minimales de { codimension positive arbitraire},  ferm\'ees et g\'eom\'etriquement distinces ?}

 \`A ce stade  du d\'eveloppement du calcul des variations de la fonctionnelle de volume, on doit reconnaitre qu'il n'y a encore aucun outil math\'ematique satisfaisant connu pour aborder cette question dans une telle g\'en\'eralit\'e. Ce probl\`eme ouvre un champs d'investigation immense et fascinant pour les d\'ecennies \`a venir.

\end{document}